

\loadbold

\define\scrO{\Cal O}
\define\Pee{{\Bbb P}}
\define\Zee{{\Bbb Z}}
\define\Cee{{\Bbb C}}
\define\Ar{{\Bbb R}}

\define\proof{\demo{Proof}}
\define\endproof{\qed\enddemo}
\define\endstatement{\endproclaim}
\define\theorem#1{\proclaim{Theorem #1}}
\define\lemma#1{\proclaim{Lemma #1}}
\define\proposition#1{\proclaim{Proposition #1}}
\define\corollary#1{\proclaim{Corollary #1}}
\define\claim#1{\proclaim{Claim #1}}

\define\section#1{\specialhead #1 \endspecialhead}
\define\ssection#1{\medskip\noindent{\bf #1}}

\documentstyle{amsppt}
\pageno=1
\topmatter
\title On complex surfaces diffeomorphic to rational surfaces
\endtitle
\author {Robert Friedman  and Zhenbo Qin}
\endauthor
\address Department of Mathematics, Columbia University, New York,
NY 10027, USA\endaddress
\email rf\@math.columbia.edu  \endemail
\address Department of Mathematics, Oklahoma State University,
Stillwater, OK 74078, USA
\endaddress
\email  qz\@math.okstate.edu \endemail
\thanks The first author was partially supported by NSF grant DMS-9203940.
The second author was partially supported by a grant from ORAU Junior Faculty
Enhancement Award Program.
\endthanks

\endtopmatter
\document

\section{Introduction}

The goal of this paper is to prove the following:

\theorem{0.1} Let $X$ be a complex surface of general type. Then $X$ is
not diffeomorphic to a rational surface.
\endstatement

Using the results from [13], we obtain the following corollary, which settles
a problem raised by Severi:

\corollary{0.2} If $X$ is a complex surface diffeomorphic to a rational
surface, then $X$ is a rational surface. Thus, up to deformation equivalence,
there is a unique complex structure on the smooth $4$-manifolds $S^2\times S^2$
and $\Pee^2 \# n\overline{\Pee}^2$.
\endstatement

In addition, as discussed in the book [15], Theorem 0.1 is the last step in
the proof of the following, which was conjectured by Van de Ven [37] (see also
[14,15]):

\corollary{0.3} If $X_1$ and $X_2$ are two diffeomorphic complex surfaces, then
$$\kappa (X_1) = \kappa (X_2),$$
where $\kappa (X_i)$ denotes the Kodaira
dimension of $X_i$.
\endstatement

The first major step in proving that every complex surface diffeomorphic to a
rational surface is  rational was Yau's theorem [40] that every complex
surface of the same homotopy type as $\Pee^2$ is biholomorphic to $\Pee^2$.
After this case, however, the problem takes on a different character: there do
exist nonrational complex surfaces with the same oriented homotopy type as
rational surfaces, and the issue is to show that they are not in fact
diffeomorphic to rational surfaces. The only known techniques for dealing with
this question involve gauge theory and  date back to Donaldson's seminal paper
[9] on the failure of the $h$-cobordism theorem in dimension 4. In this paper,
Donaldson introduced analogues of polynomial invariants for 4-manifolds $M$
with $b_2^+(M)  = 1$ and special
$SU(2)$-bundles. These invariants depend in an explicit way on a chamber
structure in the positive cone in $H^2(M; \Ar)$. Using these invariants, he
showed that a certain elliptic surface (the Dolgachev surface with multiple
fibers of multiplicities 2 and 3) was not diffeomorphic to a rational surface.
In [13], this result was generalized to cover all Dolgachev surfaces and their
blowups (the case of minimal Dolgachev surfaces was also treated in [28]) and
Donaldson's methods were also used to study self-diffeomorphisms of rational
surfaces. The only remaining complex surfaces which are homotopy equivalent
(and thus homeomorphic) to rational surfaces are then of general type, and a
single example of such surfaces, the Barlow surface, is known to exist [2]. In
1989, Kotschick [18], as well as Okonek-Van de Ven [29], using Donaldson
polynomials associated to
$SO(3)$-bundles, showed that the Barlow surface was not diffeomorphic to a
rational surface. Subsequently Pidstrigach [30] showed that no complex surface
of general type which has the same homotopy type as the Barlow surface was
diffeomorphic to a rational surface, and Kotschick [20] has outlined an
approach
to showing that no blowup of such a surface is diffeomorphic to a rational
surface. All of these approaches use $SO(3)$-invariants or $SU(2)$-invariants
for small values of the (absolute value of) the first Pontrjagin class $p_1$ of
the $SO(3)$-bundle, so that the dependence on chamber structure can be
controlled in a quite explicit way.

In [33], the second author showed that no surface $X$ of general type could be
diffeomorphic to $\Pee ^1\times \Pee ^1$ or to $\Bbb F_1$, the blowup of
$\Pee^2$ at one point. Here the main tool is the study of $SO(3)$-invariants
for large values of $-p_1$, as defined and analyzed in [19] and [21]. These
invariants also depend on a chamber structure, in a rather complicated and not
very explicitly described fashion. In [34], these methods are used to analyze
minimal surfaces $X$ of general type under certain assumptions concerning the
nonexistence of rational curves, which are always satisfied if $X$ has the same
homotopy type as $\Pee ^1\times \Pee ^1$ or $\Bbb F_1$, by a theorem of Miyaoka
on the number of rational curves of negative self-intersection on a minimal
surface of general type. The main idea of the proof is to show the following:
Let $X$ be a minimal surface of general type, and suppose that
$\{E_0, \dots, E_n\}$ is an orthogonal basis for $H^2(X; \Zee)$ with $E_0^2
=1$, $E_i^2 =-1$ for $i\geq 1$, and $[K_X] = 3E_0 -\sum _{i\geq 1}E_i$. Finally
suppose that the divisor $E_0 - E_i$ is nef for some $i\geq 1$. Then the class
$E_0 - E_i$ cannot be represented by a smoothly embedded 2-sphere. (Actually,
in [34], the proof shows that an appropriate Donaldson polynomial is not
zero whereas it must be zero if $X$ is diffeomorphic to a rational surface.
However, using [26], one can also show that if $E_0 - E_i$ is represented by a
smoothly embedded 2-sphere, then the Donaldson polynomial is zero.)  At the
same
time, building on ideas of Donaldson, Pidstrigach and Tyurin [31], using Spin
polynomial invariants, showed that no minimal surface of general type is
diffeomorphic to a rational surface.

We now discuss the contents of this paper and the general strategy for the
proof of Theorem 0.1. The  bulk of this paper is devoted to giving a new proof
of the results of Pidstrigach and Tyurin concerning minimal surfaces $X$. Here
our methods apply as well to minimal simply connected algebraic surfaces of
general type with $p_g$ arbitrary. Instead of looking at embedded 2-spheres of
self-intersection 0 as in [34], we consider those of self-intersection $-1$. We
show in fact the following in Theorem 1.10 (which includes a generalization
for blowups):

\theorem{0.4} Let $X$ be a minimal simply connected algebraic surface of
general type, and let $E\in H^2(X; \Zee)$ be a class satisfying $E^2 = -1$,
$E\cdot [K_X] = 1$. Then the class $E$ cannot be represented by a smoothly
embedded $2$-sphere.
\endstatement

In particular, if $p_g(X) = 0$, then  $X$ cannot
be diffeomorphic to a rational surface. The method of proof of Theorem 0.4 is
to show that a certain value of a Donaldson polynomial invariant for $X$ is
nonzero (Theorem 1.5), while it is a result of Kotschick that if the class $E$
is represented by a smoothly embedded 2-sphere, then the value of the Donaldson
polynomial must be zero (Proposition 1.1). In case $p_g(X) =0$, once we have
have found a polynomial invariant which distinguishes $X$ from a rational
surface, it follows in a straightforward way from the characterization of
self-diffeomorphisms of rational surfaces given in [13] that no blowup of $X$
can be diffeomorphic to a rational surface either (see Theorem 1.7). This part
of the argument could also be used with the result of Pidstrigach and Tyurin to
give a proof of Theorem 0.1.

Let us now discuss how to show that certain Donaldson polynomials do not vanish
on certain classes. The prototype for such results is the nonvanishing theorem
of Donaldson [10]: if $S$ is an algebraic surface with $p_g(S) > 0$ and
$H$ is an ample line bundle on $S$, then for all choices of $w$ and all $p \ll
0$, the
$SO(3)$-invariant $\gamma _{w,p}(H, \dots, H) \neq 0$. We give a
generalization of this result in Theorem 1.4 to certain cases where $H$ is no
longer  ample, but satisfies: $H^k$ has no base points for $k\gg 0$ and
defines a birational morphism from $X$ to a normal surface $\bar X$, and where
$p_g(X)$ is also allowed to be zero (for an appropriate choice of chamber).
Here we must assume that there is no exceptional curve $C$ such that $H\cdot
C=0$, as well as the following additional assumption concerning  the
singularities of
$\bar X$: they should be rational or minimally elliptic in the terminology of
[22].  The proof of Theorem 1.4 is a straightforward generalization of
Donaldson's original proof, together with methods developed by J\. Li in [23,
24].

Given the generalized nonvanishing theorem, the problem becomes one of
constructing divisors $M$ such that $M$ is orthogonal to a class $E$ of square
$-1$ and moreover such that $M$ is eventually base point free.
(Here we recall that a divisor $M$ is {\sl
eventually base point free\/} if the complete linear system $|kM|$ has no base
points for all
$k\gg 0$.) There are various methods for finding base point free linear systems
on an algebraic surface. For example, the well-studied method of Reider
[35] implies that, if $X$ is a minimal surface of general type and $D$ is a
nef and big divisor on
$X$, then
$M= K_X+D$  is eventually base point free. There is also a
technical generalization of this result due to Kawamata [16]. However, the
methods which we shall need are essentially elementary. The general outline of
the construction is as follows. Let $E$ be a class of square $-1$ with $K_X
\cdot E=1$. It is known that, if $E$ is the class of a smoothly embedded
$2$-sphere, then $E$ is of type $(1,1)$ [6]. Thus
$K_X+E$ is a divisor orthogonal to $E$. If $K_X+E$ is ample we are done. If
$K_X+E$ is nef but not ample, then there exist curves $D$ with $(K_X+E)\cdot
D=0$, and the intersection matrix of the set of all such curves is negative
definite. Thus we may contract the set of all such curves to obtain a normal
surface $X'$. If
$X'$ has only rational singularities, then the divisor $K_X+E$ induces a
Cartier divisor on $X'$ which is ample, by the Nakai-Moishezon criterion, and
so some multiple of $K_X+E$ is base point free. Next suppose that $X'$ has a
nonrational singular point $p$ and let $D_1, \dots, D_t$ be the
irreducible curves on $X$ mapped to $p$.
Then we give a dual form of Artin's criterion [1] for a rational
singularity, which says the following: the point $p$ is a nonrational
singularity if and only if there exist nonnegative integers $n_i$, with at
least one $n_i >0$, such that $(K_X+\sum _in_iD_i)\cdot D_j \geq 0$ for all
$j$. Moreover there is a choice of the $n_i$ such that either the inequality is
strict for every $j$ or the contraction of the
$D_j$ with $n_j \neq 0$ is a minimally elliptic singularity. In this case,
provided that $K_X$ is itself nef, it is easy to show that $K_X+\sum
_in_iD_i$ is nef and big and eventually base point free, and defines the
desired contraction. The
remaining case is when
$K_X+E$ is not nef. In this case, by considering the curves $D$ with
$(K_X+E)\cdot D<0$, it is easy to find a
$\Bbb Q$-divisor of the form
$K_X+\lambda D$, where $D$ is an irreducible curve and
$\lambda \in \Bbb Q^+$, which is nef and big and such that some
multiple is eventually base point free, and which is orthogonal to $E$. The
details are given in Section 3. These methods can also handle the case of
elliptic surfaces (the case where $\kappa (X) =1$), but of course there are
more elementary and direct arguments here which prove a more precise result.

We have included an appendix giving a proof, due to the first author, R\.
Miranda, and J\.W\. Morgan, of a result  characterizing the
canonical class of a rational surface up to isometry. This result seems to be
well-known to specialists but we were unable to find an explicit statement in
the literature. It follows from work of Eichler and Kneser on the number of
isomorphism classes of indefinite quadratic forms of rank at least 3 within a
given genus (see e\.g\. [17]) together with some calculation. However the proof
in the appendix is an elementary argument.

The methods in this paper are able to rule out the possibility of embedded
2-spheres whose associated class $E$ satisfies $E^2 = -1$, $E\cdot [K_X] = 1$.
However, in case $p_g(X) = 0$ and $b_2(X) \geq 3$, there are infinitely many
classes $E$ of square $-1$ which satisfy $|E\cdot K_X| \geq 3$. It is natural
to hope that these classes also cannot be represented by smoothly embedded
2-spheres. More generally we would like to show that the surface $X$ is
strongly minimal in the sense of [15]. Likewise, in case
$p_g(X) >0$, we have only dealt with the first case of the ``$(-1)$-curve
conjecture" (see [6]).

\medskip
\noindent
{\bf Acknowledgements:} We would like to thank
Sheldon Katz, Dieter Kotschick, and Jun Li for valuable help and stimulating
discussions.

\section{1. Statement of results and overview of the proof}

\ssection{1.1. Generalities on $\boldkey S \boldkey O\boldkey (\boldkey
3\boldkey )$-invariants}

Let $X$ be a smooth simply connected $4$-manifold with $b_2^+(X) = 1$,
and fix an $SO(3)$-bundle $P$ over $X$ with $w_2(P) = w$ and $p_1(P) = p$.
Recall that a {\sl wall of type $(w,p)$ for
$X$} is a class
$\zeta \in H^2(S; \Zee)$ such that $\zeta \equiv w \mod 2$ and
$p \leq \zeta ^2 <0$. Let
$$\Omega _X = \{x\in H^2(X; \Ar): x^2 >0\}.$$
 Let $W^\zeta = \Omega _X \cap (\zeta)^\perp$.  A
{\sl chamber of type  $(w, p)$ for $X$}  is a connected
component of the set
$$\Omega _X - \bigcup\{W^\zeta: \zeta {\text{ is a wall of type $(w,p)$}}\,
\}.$$
Let $\Cal C$ be a chamber of type $(w,p)$ for
$X$ and let
$\gamma_{w,p}(X;\Cal C)$ denote the associated Donaldson
polynomial, defined via [19] and [21]. Here $\gamma_{w,p}(X;\Cal C)$ is
only defined up to $\pm 1$, depending on the choice of an integral lift for
$w$, corresponding to a choice of orientation for the moduli space. The actual
choice of sign will not matter, since we shall only care if a certain value of
$\gamma_{w,p}(X;\Cal C)$ is nonzero. In the complex case we shall always
assume for convenience that the choice has been made so that the orientation
of the moduli space agrees with the complex orientation. Via Poincar\'e
duality, we shall view $\gamma_{w,p}(X;\Cal C)$ as a function on either
homology or cohomology classes. Given a class
$M$, we use the notation $\gamma_{w,p}(X;\Cal C) (M^d)$ for
the evaluation $\gamma_{w,p}(X;\Cal C) (M, \dots, M)$ on the class $M$ repeated
$d$ times, where
$d= -p-3$ is the expected dimension of the moduli space. We then have the
following vanishing result for
$\gamma_{w,p}(\Cal C)$, due to Kotschick [19, (6.13)]:

\proposition{1.1} Let $E \in H^2(X; \Bbb Z)$ be the cohomology class of  a
smoothly embedded $S^2$ in $X$ with $E^2 = -1$. Let $w$ be the second
Stiefel-Whitney class of $X$, or more generally any class in $H^2(X,
\Zee/2\Zee)$ such that $w\cdot E \neq 0$.  Suppose that
$M\in H_2(X; \Zee)$ satisfies $M^2 > 0$ and $M \cdot E = 0$. Then, for every
chamber $\Cal C$ of type $(w, p)$ such that the wall $W^E$ corresponding to $E$
passes through the interior of $\Cal C$,
$$\gamma_{w,p}(X;\Cal C)(M^d) = 0. \qed$$
\endproclaim

Note that if $w$ is the second Stiefel-Whitney class of $X$, then $W^E$ is
a wall of type $(w,p)$ (and so does not pass through the interior of any
chamber) if and only if $E^\perp$ is even. This case arises, for example,
if $X$ has the homotopy type of $(S^2\times S^2)\#\overline{\Pee}^2$ and $E$ is
the standard generator of $H^2(\overline{\Pee}^2; \Zee) \subseteq H^2(X; \Bbb
Z)$.

For the proof of Theorem 0.1, the result of (1.1) is sufficient. However, for
the slightly more general result of Theorem 1.10, we will also need the
following variant of (1.1):

\theorem{1.2} Let $E \in H^2(X; \Bbb Z)$ be the cohomology class of  a
smoothly embedded $S^2$ in $X$ with $E^2 = -1$. Let $w$ be a class in
$H^2(X, \Zee/2\Zee)$ such that $w\cdot E \neq 0$.  Suppose that
$M\in H_2(X; \Zee)$ satisfies $M^2 > 0$ and $M \cdot E = 0$. Then, for every
chamber $\Cal C$ of type $(w, p)$ containing $M$ in its closure,
$$\gamma_{w,p}(X;\Cal C)(M^d) = 0,$$
unless $p=-5$ and $w$ is the \rom{mod} $2$ reduction of $E$.
Thus, except in this last case, $\gamma_{w,p}(X;\Cal C)$ is divisible by $E$.
\endproclaim
\proof If $W^E$ is not a wall of type $(w,p)$ we are done by (1.1). Otherwise,
$E$ defines a wall of type $(w,p)$ containing $M$. Next let us assume that
$E^\perp \cap \overline{\Cal C}$ is a codimension one face of the
closure $\overline{\Cal C}$ of $\Cal C$. We have an induced
decomposition of $X$:
$$X = X_0 \# \overline{\Pee}^2.$$
Identify $H_2(X_0; \Bbb Z)$ with the subspace $E^\perp$ of $H_2(X; \Bbb Z)$,
and let $\overline{\Cal C}_0 = E^\perp \cap \overline{\Cal C}$. Then
$\overline{\Cal C}_0$ is the closure of some chamber ${\Cal C}_0$
of type $(w - e, p + 1)$ on $X_0$, where $e$ is the mod 2 reduction of $E$.
Choose a generic Riemannian metric $g_0$ on $X_0$ such that the cohomology
class $\omega _0$ of the self-dual harmonic 2-form associated to
$g_0$ lies in the interior of
$\Cal C_0$. By the results in [39], there is a family of metrics $h_t$ on the
connected sum $X_0\# \overline{\Cee P}^2$ which converge in an appropriate
sense to $g_0\amalg g_1$, where $g_1$ is the Fubini-Study metric on
$\overline{\Pee}^2$, and such that the cohomology classes of the self-dual
harmonic 2-forms associated to $h_t$ lie in the interior of
$\Cal C$ and converge to $\omega _0$.

Standard gluing and compactness arguments (see for example [15],
Appendix to Chapter 6) and dimension counts show that the restriction
of the invariant
$\gamma_{w,p}(X;\Cal C)$ to $H_2(X_0; \Bbb Z)$ vanishes.

Consider now the general case where $W^E$ is a wall of type $(w,p)$ and the
closure of $\Cal C$ contains $M$ but where
$W^E\cap \overline{\Cal C}$ is not necessarily a codimension one face of
$\overline{\Cal C}$. Since
$W^E$ is a wall of type $(w,p)$ and $M\in E^\perp$, there exists a chamber
$\Cal C'$ of type $(w, p)$ whose closure contains $M$ such that $W^E$ is a
codimension one face of $\overline{\Cal C}'$. By the previous argument,
$\gamma_{w,p}(X;\Cal C')(M^d) =0$ and so it
will suffice to show that
$$\gamma_{w,p}(X;\Cal C)(M^d) = \gamma_{w,p}(X;\Cal C')(M^d).$$
Note that $\Cal C$ and $\Cal C'$ are separated by  finitely many walls of type
$(w, p)$ all of which contain the class $M$.  Thus, we have a sequence of
chambers of type $(w, p)$:
$$\Cal C = \Cal C_1, \Cal C_2, \dots, \Cal C_{k - 1},
\Cal C_k = \Cal C'$$
such that for each $i$, $\Cal C_{i - 1}$ and $\Cal C_i$ are
separated by a single wall $W_i=W^{\zeta_i}$ of type $(w, p)$ which contains
$M$. Since $W_i$ contains $M$, $M \cdot \zeta_i = 0$. By [19, (3.2)(3)] (see
also [21]), the difference $\gamma_{w,p}(X;\Cal C_{i - 1}) -
\gamma_{w,p}(X;\Cal C_i)$ is divisible by the class $\zeta_i$ except in the
case where $p=-5$ and $w$ is the mod 2 reduction of $E$. It follows that,
except in this last case, for each $i$,
$$\gamma_{w,p}(X;\Cal C_{i - 1})(M^d) = \gamma_{w,p}(X;\Cal C_i)(M^d).$$
Hence $\gamma_{w,p}(X;\Cal C)(M^d) =
\gamma_{w,p}(X;\Cal C')(M^d)=0$.
\endproof

We shall also need the following ``easy" blowup formula:

\lemma{1.3} Let $X\#\overline{\Pee}^2 $ be a blowup of $X$, and identify
$H_2(X; \Bbb Z)$ with a subspace of
$H_2(X\#\overline{\Pee}^2; \Bbb Z)$ in the natural way.
Given $w\in H^2(X; \Zee/2\Zee)$, let $\tilde \Cal C$ be a chamber of
type $(w, p)$ for $X\#\overline{\Pee}^2$ containing
the chamber $\Cal C$ in its closure. Then
$$\gamma _{w,p}(X\#\overline{\Pee}^2;\tilde\Cal C)
|H_2(X; \Bbb Z) = \pm\gamma_{w,p}(X;\Cal C).$$
\endproclaim
\proof Choose a generic Riemannian metric $g$ on
$X$ such that the cohomology class $\omega$ of the self-dual harmonic 2-form
associated to $g$ lies in the interior of
$\Cal C$. We again use the results in [39] to choose a family of metrics $h_t$
on the connected sum $X\# \overline{\Pee}^2$ which converge in an appropriate
sense to $g\amalg g'$, where $g'$ is the Fubini-Study metric on
$\overline{\Pee}^2$, and such that the cohomology classes of the self-dual
harmonic 2-forms associated to $h_t$ lie in the interior of
$\tilde\Cal C$ and converge to $\omega$. Standard gluing and compactness
arguments (see e\.g\. [15], Chapter 6, proof of Theorem 6.2(i)) now show that
the restriction of $\gamma _{w,p}(X\#\overline{\Pee}^2;\tilde\Cal C)$ to
$H_2(X; \Bbb Z)$ (with the appropriate orientation conventions) is just
$\gamma_{w,p}(X;\Cal C)$.
\endproof

\ssection{1.2. The case of a minimal $\boldkey X$}

In this subsection we shall outline the results to be proved concerning minimal
surfaces of general type. One basic tool is a nonvanishing theorem for certain
values of the Donaldson polynomial:

\theorem{1.4}  Let $X$ be a simply connected algebraic surface with $p_g(X)
=0$,
and let $M$ be a nef and big divisor on $X$ which is eventually
base point free. Denote by $\varphi\: X \to \bar X$ the birational morphism
defined by $|kM|$ for $k\gg 0$, so that $\bar X$ is a normal projective
surface.  Suppose that $\bar X$ has only rational or minimally elliptic
singularities, and that $\varphi$ does not contract any exceptional curves to
points. Let
$w\in H^2(X;\Zee/2\Zee)$ be the \rom{mod} $2$ reduction of the class $[K_X]$.
Then there exists a constant $A$ depending only on $X$ and $M$ with the
following property:
For all integers $p\leq A$, let
$\Cal C$ be a chamber of type $(w,p)$ containing $M$ in its closure and suppose
that $\Cal C$ has nonempty intersection with the ample cone of $X$. Set $d
= -p-3$. Then
$$\gamma_{w,p}(X;\Cal C)(M^d) >0.$$
\endstatement

We shall prove Theorem 1.4 in Section 2, where we shall also recall the salient
properties of rational and minimally elliptic singularities. The proof also
works in  the case where
$p_g(X)>0$, in which case
$\gamma_{w,p}(X)$ does not depend on the choice of a chamber.

We can now state the main result concerning minimal surfaces, which we shall
prove in Section 3:

\theorem{1.5} Let $X$ be a minimal simply connected algebraic surface of
general type, and let
$E\in H^2(X; \Zee)$ be a $(1,1)$-class satisfying $E^2=-1$, $E\cdot K_X = 1$.
Let $w$ be the \rom{mod} $2$ reduction of $[K_X]$. Then there exist:
\roster
\item"{(i)}"  an integer $p$
and \rom(in case $p_g(X)=0$\rom) a chamber $\Cal C$ of type $(w,p)$ and
\item"{(ii)}" a $(1,1)$-class $M\in H^2(X; \Zee)$
\endroster
such that
$M\cdot E=0$ and $\gamma _{w,p}(X)(M^d) \neq 0$
\rom(or, in case $p_g(X)=0$,
$\gamma _{w,p}(X; \Cal C)(M^d) \neq 0$\rom).
\endstatement

The method of proof of (1.5) will be the following: we will show that there
exists an orientation preserving self-diffeomorphism $\psi$ of
$X$ with $\psi ^*[K_X] = [K_X]$ and  a nef and big divisor $M$ on $X$ such
that:
\roster
\item"{(i)}" $M\cdot \psi ^*E = 0$.
\item"{(ii)}" $M$ is eventually base point free,
and the corresponding contraction $\varphi\: X \to \bar X$ maps $X$
birationally onto a normal surface $\bar X$ whose only singularities are either
rational or minimally elliptic.
\endroster
Using the naturality of $\gamma _{w,p}(X;\Cal C)$, it suffices to prove (1.5)
after replacing $E$ by $\psi ^*E$. In this case, by Theorem 1.4 with
$w= [K_X]$,
$\gamma _{w,p}(X;\Cal C)(M^d)
\neq 0$ for all
$p \ll 0$.

\corollary{1.6} Let $X$ be a simply connected minimal surface of general type
with $p_g(X)=0$.  Then there exist
\roster
\item"{(i)}" a class $w\in H^2(X; \Zee/2\Zee)$;
\item"{(ii)}" an integer $p\in \Zee$;
\item"{(iii)}" a chamber $\Cal C$ for $X$ of type $(w,p)$, and
\item"{(iv)}" a homotopy equivalence $\alpha: X \rightarrow Y$, where $Y$ is
either the blowup of ${\Pee}^2$  at $n$ distinct points or $Y=\Pee ^1\times
\Pee ^1$,
\endroster
such that, for $w'= (\alpha ^*)^{-1}(w)$ and $\Cal C' = (\alpha ^*)^{-1}(\Cal
C)$,
$$\alpha ^*\gamma _{w',p}(Y;\Cal C') \neq
\pm\gamma_{w, p}(X; \Cal C).$$
\endproclaim
\proof  If $X$
is homotopy equivalent to $\Pee ^1\times \Pee ^1$ then the theorem
follows from [33]. Otherwise $X$ is oriented homotopy equivalent to $\Pee ^2\#
n\overline{\Pee}^2$, for $1\leq n\leq 8$, and we claim that there exists a
homotopy equivalence $\alpha
\: X \to Y$ such that $\alpha ^*[K_Y] = -[K_X]$. Indeed, every
integral isometry $H^2(Y; \Zee) \to H^2(X; \Zee)$ is realized by an oriented
homotopy equivalence. Thus it suffices to show that every two characteristic
elements of
$H^2(Y; \Zee)$ of square $9-n$ are conjugate under the isometry group, which
follows from the appendix to this paper.  Choosing
such a homotopy equivalence $\alpha$, let $e$ be the class of an exceptional
curve in $Y$ and let $E = \alpha ^*e$. Then $E^2 =-1$ and $E\cdot [K_X] = 1$.
We may now apply Theorem 1.5 to the class
$E$, noting that
$E$ is a $(1,1)$ class since $p_g(X)=0$. Let $M$ and $\Cal C$ be a divisor
and a chamber which satisfy the conclusions of Theorem 1.5 and let
$m = (\alpha ^*)^{-1} M$. If $w$ is the mod 2 reduction of $[K_X]$, then
$w'$ is the mod two reduction of $[K_Y]$, so that $w'$ is also
characteristic. Now
$e$ is the class of a smoothly embedded 2-sphere in $Y$ since it is the class
of an exceptional curve. Moreover
$m\cdot e =0$. By
Theorem 1.2,
$\gamma_{w', p}(Y;\Cal C')(m^d) = 0$ since $e$ is represented
by a smoothly embedded 2-sphere and $w'$ is characteristic. But $\gamma
_{w,p}(X;\Cal C)(M^d)
\neq  0$ by Theorem 1.5 and the choice of $M$. Thus
$\alpha ^*\gamma _{w',p}(Y;\Cal C') \neq
\pm\gamma_{w, p}(X; \Cal C)$.
\endproof

Using the result of Wall [38] that every homotopy self-equivalence from $Y$ to
itself is realized by a diffeomorphism, the proof above shows that the
conclusions of the corollary hold for {\it every\/} homotopy equivalence
$\alpha: X \rightarrow Y$.

\ssection{1.3. Reduction to the minimal case}

We begin by recalling some terminology and results from [13].
A {\sl good generic rational surface} $Y$ is a rational surface
such that $K_Y = -C$ where $C$ is a smooth curve,
and such that there does not exist a smooth rational curve on $Y$
with self-intersection $-2$. Every rational surface is
diffeomorphic to a good generic rational surface.

\theorem{1.7} Let $X$ be a minimal surface of general type and
let $\tilde X\to X$ be a blowup of $X$ at $r$ distinct points.
Let $E_1', \dots , E_r'$ be the homology classes of
the exceptional curves on $\tilde X$.
Let $\psi _0\: \tilde X \to \tilde Y$ be a diffeomorphism,
where $\tilde Y$ is a good generic rational surface.
Then there exist a diffeomorphism $\psi \: \tilde X \to \tilde Y$
and a good generic rational surface $Y$ with the following properties:
\roster
\item"{(i)}" The surface $\tilde Y$ is the blowup of $Y$ at
$r$ distinct points.
\item"{(ii)}" If $e_1, \dots, e_r$ are the classes of
the exceptional curves in $H^2(\tilde Y; \Bbb Z)$ for the blowup $\tilde Y \to
Y$,  then possibly after renumbering $\psi^*(e_i) = E_i'$ for all $i$.
\item"{(iii)}" Identifying $H^2(X)$ with a subgroup of $H^2(\tilde X)$
and $H^2(Y)$ with a subgroup of $H^2(\tilde Y)$ in the obvious way,
we have $\psi^*(H^2(Y)) = H^2(X)$.
\endroster
Moreover, for every choice of an isometry $\tau$ from $H^2(Y)$ to $H^2(X)$,
there exists a choice of a diffeomorphism $\psi$ satisfying \rom{(i)--(iii)}
above and such that $\psi ^*|H^2(Y) = \tau$.
\endstatement
\proof Let $e_i' \in H^2(\tilde Y; \Zee)$ satisfy $\psi_0^*(e_i') = E_i'$.
Thus the Poincar\'e dual of $e_i'$ is represented by
a smoothly embedded 2-sphere in $\tilde Y$.
It follows that reflection $r_{e_i'}$ in $e_i'$ is realized
by an orientation-preserving self-diffeomorphism of $\tilde Y$.
To see what this says about $e_i'$, we shall recall the following
terminology from [13].

Let $\bold H(\tilde Y)$ be the set
$\{\, x\in H^2(\tilde Y; \Ar)\mid x^2 =1\,\}$ and
let $\Cal K(\tilde Y)\subset H^2(\tilde Y; \Ar)$ be intersection of
the closure of the K\"ahler cone of $\tilde Y$ with $\bold H(\tilde Y)$.
Then $\Cal K(\tilde Y)$ is a convex subset of $\bold H(\tilde Y)$
whose walls consist of the classes of exceptional curves on $\tilde Y$
together with $[-K_{\tilde Y}]$ if $b_2^-(\tilde Y)
\geq 10$, which is confusingly called the {\sl exceptional wall\/} of $\Cal
K(\tilde Y)$.  Let
$\Cal R$ be the group generated by the reflections in the walls  of $\Cal
K(\tilde Y)$ defined by exceptional classes  and define the super $P$-cell
$\bold S = \bold S(P)$ by
$$\bold S = \bigcup _{\gamma \in \Cal R}\gamma \cdot \Cal K(\tilde Y).$$
By Theorem 10A on p\. 355 of [13], for an integral isometry
$\varphi$ of $H^2(\tilde Y; \Ar)$, there exists a diffeomorphism of
$\tilde Y$ inducing $\varphi$ if and only if
$\varphi (\bold S) = \pm \bold S$. (Here, if $b_2^-(\tilde Y) \leq 9$,
$\bold S = \bold H$ and the result reduces to a result
of C\.T\.C\. Wall [38].) Note that $\bold H(\tilde Y)$ has two connected
components, and reflection $r_e$ in a class $e$ of square $-1$ preserves the
set of connected components. Thus if $r_e(\bold S) = \pm \bold S$, then
necessarily
$r_e(\bold S) = \bold S$.

Next we have the following purely algebraic lemma:

\lemma{1.8} Let $e$ be a class of square $-1$ in $H^2(\tilde Y; \Zee)$ such
that the reflection $r_e$ satisfies $r_e(\bold S) = \bold S$. Then there is an
isometry $\varphi$ of $H^2(\tilde Y; \Zee)$ preserving $\bold S$ which sends
$e$ to the class of an exceptional curve.
\endstatement
\proof We first claim that, if $W$ is the wall corresponding to $e$,
then $W$ meets the interior of $\bold S$. Indeed,
the interior $\operatorname{int}\bold S$ of $\bold S$ is connected, by
Corollary 5.5 of [13] p\. 340.  If $W$ does not meet $\operatorname{int} \bold
S$, then the sets
$$\{\, x\in \operatorname{int}\bold S\mid e\cdot x > 0\,\}$$
and
$$\{\, x\in \operatorname{int}\bold S\mid e \cdot x < 0\,\}$$
are disjoint open sets covering $\operatorname{int}\bold S$
which are exchanged under the reflection $r_e$.
Since at least one is nonempty, they are both nonempty,
contradicting the fact that $\operatorname{int} \bold S$ is connected.
Thus $W$ must meet $\operatorname{int}\bold S$.

Now let $C$ be a chamber for the walls of square $-1$ which has $W$ as a wall.
It follows from Lemma 5.3(b) on p\. 339 of [13] that $C\cap \bold S$ is a
$P$-cell $P$ and that
$W$ defines a wall of $P$ which is not the exceptional wall.
By Lemma 5.3(e) of [13],
$\bold S$ is the unique super $P$-cell containing $P$, and the reflection group
generated by the elements of square $-1$ defining the walls of $P$ acts simply
transitively on the $P$-cells in $\bold S$. There is thus an element $\varphi$
in this reflection group which preserves $\bold S$ and sends $P$ to $\Cal
K(\tilde Y)$ and $W$ to a  wall of $\Cal K(\tilde Y)$ which is not
an exceptional wall.  It follows that $\varphi (e)$ is the class of an
exceptional curve on $\tilde Y$.
\endproof

Returning to the proof of Theorem 1.7, apply the previous lemma to the
reflection in $e_r'$. There is thus an isometry $\varphi$ preserving $\bold S$
such that $\varphi(e_r') = e_r$, where $e_r$ is the class of an exceptional
curve on $\tilde Y$.  Moreover
$\varphi$ is realized by a diffeomorphism.

Thus after composing with the diffeomorphism inducing $\varphi$,
we can assume that $e_r' = e_r$, or equivalently that $\psi _0^*e_r = [E_r']$.
Let $\tilde Y \to \tilde Y_1$ be  the blowing down of the exceptional curve
whose class is $e_r$. Then
$\tilde Y_1$ is again a good generic surface  by [13] p\. 312 Lemma 2.3. Since
$e_1', \dots, e_{r-1}'$ are orthogonal to
$e_r$,  they lie in the subset $H^2(\tilde Y_1)$ of $H^2(\tilde Y)$.
For $i \neq r$, the reflection in $e_i'$ preserves $W\cap \bold S$,
where $W= (e_r)^\perp$. Now $W$ is just $H^2(\tilde Y_1)$ and
$\Cal K(\tilde Y) \cap H^2(\tilde Y_1) = \Cal K(\tilde Y_1)$
by [13] p\. 331 Proposition 3.5. The next lemma relates the corresponding super
$P$-cells:

\lemma{1.9} $W\cap \bold S$ is the super $P$-cell
$\bold S_1$ for $\tilde Y_1$ containing $\Cal K(\tilde Y_1)$.
\endstatement
\proof Trivially $\bold S_1 \subseteq W\cap \bold S$,
and both sets are convex subsets with nonempty interiors.
If they are not equal, then there is a $P$-cell $P' \subset \bold S_1$
and an exceptional wall of $P'$ which passes through the interior
of $\bold S\cap W$. If $\kappa (P')$ is the exceptional wall meeting
$\bold S\cap W$, then, by [13] p\. 335 Lemma 4.6, $\kappa (P') - e_r$ is
an exceptional wall of $P$ for a well-defined $P$-cell in $\bold S$,
and $\kappa (P') - e_r$ must pass through the interior of $\bold S$.
This is a contradiction. Hence $\bold S\cap W = \bold S_1$ is
a super $P$-cell of $\tilde Y_1$, and we have seen that it contains $\Cal
K(\tilde Y_1)$.
\endproof

Returning to the proof of Theorem 1.7, reflection in $e_{r-1}'$
preserves $\bold S_1$. Applying Lemma 1.8,  there is a diffeomorphism of
$\tilde Y_1$  which sends $e_{r-1}'$ to the class of an exceptional curve
$e_{r-1}$.  Of course, there is an induced diffeomorphism of $\tilde Y$
which fixes $e_r$. Now we can clearly proceed by induction on $r$.

The above shows that after replacing $\psi _0$ by
a diffeomorphism $\psi$ we can find $Y$ as above
so that (i) and (ii) of the statement of Theorem 3 hold.
Clearly $\psi^*(H^2(Y)) = H^2(X)$. By the theorem of C\.T\.C\. Wall mentioned
above, there is a diffeomorphism of $Y$  realizing every integral isometry of
$H^2(Y)$.  So after further modifying by a diffeomorphism of $Y$,
which extends to a diffeomorphism of $\tilde Y$ fixing
the classes of the exceptional curves, we can assume that
the diffeomorphism $\psi$ restricts to $\tau$ for
any given isometry from $H^2(Y)$ to $H^2(X)$.
\endproof

We can now give a proof of Theorem 0.1:

\theorem{0.1} No complex surface of general type is diffeomorphic to a rational
surface.
\endstatement
\proof Suppose that $X$ is a minimal surface of general type and that
$\rho \: \tilde X \to X$ is a blowup of $X$ diffeomorphic to a rational
surface.
We may assume that $\tilde X$ is diffeomorphic via $\psi$
to a good generic rational surface ${\tilde Y}$,
and that $\rho '\: {\tilde Y} \to Y$ is a blow up of ${\tilde Y}$
such that $Y$ and $\psi$ satisfy (i)--(iii) of Theorem 1.7. Choose
$w, p, \alpha, \Cal C$ for
$X$ such that the conclusions of  Corollary 1.6 hold, with $\Cal C'$ the
corresponding chamber on $Y$, and let
$\tilde
\Cal C'$ be any chamber  for ${\tilde Y}$ containing $\Cal C'$ in its closure.
Then $\psi ^*\tilde \Cal C' = \tilde \Cal C$ is a chamber on $\tilde X$
containing $\Cal C$ in its closure. Using the last sentence of Theorem 1.7, we
may  assume that
$\psi ^*|H^2(Y) = \alpha^*$. Thus $\psi ^* (\rho ') ^* = \rho ^*\alpha ^*$. By
the functorial properties of  Donaldson polynomials, and viewing $H^2(X;
\Zee/2\Zee)$ as a subset of $H^2(\tilde X; \Zee/2\Zee)$, and similarly for
$\tilde Y$, we have
$$\psi ^*\gamma _{w', p}({\tilde Y}, \tilde \Cal C')
= \pm\gamma _{\psi ^*w', p}(\tilde X, \tilde \Cal C)=
\pm\gamma _{w, p}(\tilde X, \tilde \Cal C).$$
Restricting each side to $\psi ^*H_2(Y) = H_2(X)$,
we obtain by repeated application of Lemma 1.3 that
$$\alpha ^*\gamma _{w',p}(Y;\Cal C') = \pm\gamma _{w, p}(X;\Cal C).$$
But this contradicts Corollary 1.6.
\endproof

Using Theorem 1.5, we have the following generalization of Theorem 0.4 in the
introduction to the case of nonminimal algebraic surfaces:

\theorem{1.10} Let $X$ be a minimal simply connected surface of general type,
and let $E\in H^2(X; \Zee)$ satisfy $E^2 = -1$ and
$E\cdot K_X = 1$. Let $\tilde X$ be a blowup of $X$. Then, viewing $H^2(X;
\Zee)$ as a subset of $H^2(\tilde X; \Zee)$, the class $E$ is not represented
by a smoothly embedded $2$-sphere in $\tilde X$.
\endstatement
\proof Suppose instead that $E$ is represented by a smoothly embedded
$2$-sphere. If
$p_g(X) >0$, then it follows from the results of [6] that $E$ is a
$(1,1)$-class, i\.e\. $E$ lies in the image of
$\operatorname{Pic}X$ inside $H^2(X; \Zee)$. Of course, this is automatically
true if $p_g(X) = 0$. Next assume that $p_g(X) = 0$. By  Theorem 1.5,
there exists a
$w\in H^2(X; \Zee/2\Zee)$, an integer $p$, and a chamber
$\Cal C$ of type $(w,p)$, such that $\gamma _{w,p}(X;\Cal C)(M^d)\neq
0$, where $M$ is a class in the closure of $\Cal C$ and $M\cdot E = 0$.
Consider the Donaldson polynomial
$\gamma_{w,p}(\tilde X;\tilde \Cal C)$, where we view $w$ as an element
of $H^2(\tilde X;\Zee/2\Zee)$ in the natural way and $\tilde
\Cal C$ is a chamber of type $(w,p)$ on $\tilde X$ containing $\Cal C$ in its
closure. Then $\tilde
\Cal C$ also contains $M$ in its closure. Thus, by Theorem 1.2,
$\gamma_{w,p}(\tilde X;\tilde \Cal C)(M^d) = 0$. On the other hand, by Lemma
1.3, $\gamma_{w,p}(\tilde X;\tilde
\Cal C)(M^d) = \pm \gamma _{w,p}(X;\Cal C)(M^d)\neq
0$. This is a contradiction. The case where $p_g(X) > 0$ is similar.
\endproof

We also have the following corollary, which works under the assumptions of
Theorem 1.10 for surfaces with $p_g>0$:

\corollary{1.11} Let $X$ be a  simply connected surface of general type with
$p_g(X) >0$, not necessarily minimal, and let $E\in H^2(X; \Zee)$ satisfy
$E^2 = -1$ and $E\cdot K_X = -1$. Suppose that $E$ is represented by a smoothly
embedded $2$-sphere. Then $E$ is the cohomology class
associated to an exceptional curve.
\endstatement
\proof Using [15] and [6], we see that if $E$ is not the cohomology class
associated to an exceptional curve, then $E\in H^2(X_{\text{min}}; \Zee)$,
where $X_{\text{min}}$ is the minimal model of $X$ and we have the natural
inclusion $H^2(X_{\text{min}}; \Zee) \subseteq H^2(X; \Zee)$. We may then apply
Theorem 1.10 to conclude that $-E$ cannot be represented by a smoothly embedded
$2$-sphere, and thus that $E$ cannot be so represented, a contradiction.
\endproof

\section{2. A generalized nonvanishing theorem}

\ssection{2.1. Statement of the theorem and the first part of the proof}

In this section, we shall prove Theorem 1.4. We first recall its statement:

\theorem{1.4} Let $X$ be a simply connected algebraic surface with $p_g(X) =0$,
and let $M$ be a nef and big divisor on $X$ which is eventually
base point free. Denote by $\varphi\: X \to \bar X$ the birational morphism
defined by $|kM|$ for $k\gg 0$, so that $\bar X$ is a normal projective
surface.
Suppose that $\bar X$ has only rational or minimally elliptic singularities,
and that $\varphi$ does not contract any exceptional curves to points. Let
$w\in H^2(X;\Zee/2\Zee)$ be the \rom{mod} $2$ reduction of the class $[K_X]$.
Then there exists a constant $A$ depending only on $X$ and $M$ with the
following property:
For all integers $p\leq A$, let
$\Cal C$ be a chamber of type $(w,p)$ containing $M$ in its closure and suppose
that $\Cal C$ has nonempty intersection with the ample cone of $X$. Set $d
= -p-3$. Then
$$\gamma_{w,p}(X;\Cal C)(M^d) >0.$$
A similar conclusion holds if $p_g(X) >0$.
\endstatement
\proof We begin by fixing some notation.  For
$L$ an ample line bundle on
$X$, given  a divisor $D$ on $X$ and an
integer $c$, let $\frak M_L(D, c)$ denote the moduli space of isomorphism
classes of
$L$-stable rank two holomorphic vector bundles on $X$ with $c_1(V) = D$ and
$c_2(V) = c$. Let $w$ be the mod 2 reduction of $D$ and let $p= D^2 - 4c$. Then
we also denote $\frak M_L(D, c)$ by
$\frak M_L(w, p)$, the moduli space of equivalence classes of $L$-stable
rank two holomorphic vector bundles on $X$ corresponding to the choice of
$(w,p)$. Here we recall that two vector bundles $V$ and $V'$ are {\sl
equivalent\/} if there exists a holomorphic line bundle $F$ such that $V'
= V\otimes F$. The invariants $w$ and $p$ only depend on the equivalence
class of $V$. Let
$\overline{\frak M_L(w, p)}$ denote the Gieseker compactification of
$\frak M_L(w, p)$, i\.e\. the
Gieseker compactification $\overline{\frak M_L(D, c)}$ of $\frak M_L(D, c)$.
Thus $\overline{\frak M_L(w, p)}$ is a projective variety.

We now fix a compact neighborhood $\Cal N$ of $M$ inside the positive cone
$\Omega _X$ of
$X$. Note that, since $M$ is nef, such a neighborhood has nontrivial
intersection with the ample cone of $X$. Using a straightforward extension of
the theorem of Donaldson [10] on the dimension of the moduli space (see e\.g\.
[12] Chapter 8, [32], [24]), there exist constants $A$ and $A'$ such that, for
all  ample line bundles $L$ such that $c_1(L) \in \Cal N$,  the following
holds:
\roster
\item If $p \leq A$, then the moduli space $\overline{\frak M_L(w, p)}$ is
good, in other words it is generically reduced of the correct dimension $-p-3$;
\item $\frak M_L(w, p)$ is a dense open subset of $\overline{\frak M_L(w, p)}$
and the generic point of $\overline{\frak M_L(w, p)}- \frak M_L(w, p)$
correspond to a torsion free sheaf $V$ such that the length of
$V\spcheck{}\spcheck/V$ is one and such that the support of
$V\spcheck{}\spcheck/V$ is a generic point of $X$;
\item For all $p' \geq A$, the dimension $\dim \frak M_L(w, p') \leq A'$.
\endroster

We shall need to make one more assumption on the integer $p$. Let $\varphi\: X
\to \bar X$ be the contraction morphism associated to $M$. For each connected
component $E$ of the set of exceptional fibers of $\varphi$, fix a (possibly
nonreduced) curve $Z$ on $X$ whose support is exactly $E$. In practice we shall
always take $Z$ to be the fundamental cycle of the singularity, to be defined
in Subsection 2.3 below. A slight generalization ([12], Chapter 8) of
Donaldson's theorem on the dimension of the moduli space then shows the
following: after possibly modifying the constant $A$,
\roster
\item"(4)" The generic
$V\in \frak M_L(w,p)$ satisfies: the natural map
$$H^1(X; \operatorname{ad}V) \to H^1(Z; \operatorname{ad}V|Z)$$
is surjective. In other words, the local universal deformation of $V$ is versal
when viewed as a deformation of $V|Z$ (keeping the determinant fixed).
\endroster

We now assume that $p\leq A$.
Let $L$ be
an ample line bundle which is not separated from $M$ by any wall of type
$(w,p)$ (or equivalently of type $(D, c)$), and moreover does not lie on
any wall of type $(w,p)$.  Thus by assumption, none of the points of
$\overline{\frak M_L(D, c)}$ corresponds to a strictly semistable sheaf.
Let $C\subset
X$ be a smooth curve  of genus
$g$. Suppose that $C\cdot D = 2a$ is even. Choosing a line bundle $\theta$ of
degree $g-1-a$ on
$C$, we can form the determinant line bundle $\Cal L(C, \theta)$ on the moduli
functor associated to torsion free sheaves corresponding to the values $w$ and
$p$ ([15], Chapter 5). Using Proposition 1.7 in [23], this line bundle descends
to a line bundle on
$\overline{\frak M_L(w, p)}$, which we shall continue to denote by $\Cal L(C,
\theta)$. Moreover, by the method of proof of Theorem 2 of [23], the line
bundle $\Cal L(C, \theta)$ depends only on the linear equivalence class of $C$,
in the sense that if $C$ and $C'$ are linearly equivalent and $\theta '$ is a
line bundle of degree $g-1-a$ on $C'$, then $\Cal L(C, \theta) \cong \Cal L(C',
\theta')$.

Next we shall use the following result, whose proof is deferred to the next
subsection:

\lemma{2.1} In the above notation, if $k\gg 0$ and $C\in |kM|$ is a smooth
curve, then, for all $N\gg 0$, the linear system associated  to $\Cal L(C,
\theta)^N$ has no base points and defines a generically finite morphism from
$\overline{\frak M_L(w, p)}$ to its image. In particular, if $d = \dim
\overline{\frak M_L(w, p)}$, then
$$c_1(\Cal L(C, \theta))^d > 0.$$
\endstatement

It follows by applying an easy adaptation of Theorem 6 in [23] or the results
of [25] to the case $p_g(X) =0$ that, since the spaces $\frak M_L(w,p')$ have
the expected dimension for an appropriate range of $p'\geq p$,
$c_1(\Cal L(C, \theta))^d$ is exactly the value  $k^d\gamma_{w,p}(X;\Cal
C)(M^d)$. Thus we have proved
Theorem 1.4, modulo the proof of Lemma 2.1. This proof will be given below.
\enddemo

\ssection{2.2. A generalization of a result of Bogomolov}

We keep the notation of the preceding subsection. Thus $M$ is a nef and big
divisor  such  that the complete linear system $|k M|$ is base point free
whenever $k \gg 0$. Throughout, we shall further assume that $M$ is divisible
by $2$ in $\operatorname{Pic}X$. Moreover
$w$ and $p$ are now fixed and $L$ is an ample line bundle such that $c_1(L) \in
\Cal N$ is not separated from $M$ by a wall of type $(w,p)$ and moreover that
$c_1(L)$ does not lie on a wall of type $(w,p)$. In particular the determinant
line bundle $\Cal L(C,
\theta)$ is defined for all smooth $C$ in $|kM|$ for all $k\gg 0$.

We then have the following generalization of a restriction theorem due to
Bogomolov [4]:

\lemma{2.2} With the above notation, there exists a constant $k_0$ depending
only on $w$, $p$, $M$, and $L$, such that for all
$k\geq k_0$ and all smooth curves
$C\in |kM|$, the following holds: for all $c' \leq c$ and $V \in \frak M_L(D,
c')$, either $V|C$ is semistable or there exists a divisor $G$
on $X$, a zero-dimensional subscheme
$\Cal Z$ and an exact sequence
$$0 \to \scrO_X(G) \to V \to \scrO_X(D-G) \otimes I_{\Cal Z} \to 0,$$
where $2G-D$ defines a wall of type $(w,p)$ containing $M$ and $C \cap
\operatorname{Supp}\Cal Z \neq \emptyset$.
\endstatement
\proof The proof follows closely the original proof of Bogomolov's theorem
[4] or [15] Section 5.2.
Choose $k_0 \geq -p$ and assume also that there exists a smooth curve $C$ in
$|kM|$ for all $k\geq k_0$. Suppose that $V|C$ is not semistable. Then there
exists a surjection $V|C \to F$, where $F$ is a line bundle on $C$ with $\deg
F=f< (D\cdot C)/2$. Let $W$ be the kernel of the induced surjection $V\to F$.
Thus
$W$ is locally free and there is an exact sequence
$$0 \to W \to V \to F \to 0.$$
A calculation gives
$$\align
p_1(\operatorname{ad}W) &= p_1(\operatorname{ad}V) + 2D\cdot C + (C)^2 -
4f\\
&> p + k^2(M)^2 \geq p + p^2 \geq 0.
\endalign$$
By Bogomolov's inequality, $W$ is unstable with respect to every ample line
bundle on
$X$. Thus there exists a divisor $G_0$ and an injection $\scrO_X(G_0) \to W$
(which we may assume to have torsion free cokernel) such that
$2(L\cdot G_0) > L\cdot (D-C)$, i\.e\. $L \cdot (2G_0 - D +C) > 0$. By
hypothesis there is an exact sequence
$$0 \to \scrO_X(G_0) \to W \to \scrO_X(-G_0+D-C) \otimes I_{\Cal Z_0}\to 0.$$
Thus
$$0< p_1(\operatorname{ad}W) = (2G_0 - D +C)^2 - 4\ell (\Cal Z_0) \leq
(2G_0 - D +C)^2.$$
It follows that $(2G_0 - D +C)^2>0$.
As  $L \cdot (2G_0 - D +C) > 0$ and $(2G_0 - D +C)^2 > 0$, $M \cdot (2G_0 - D
+C) \geq 0$ as well, i\.e\. $-(M \cdot (2G_0 - D))\leq k(M)^2$. On the other
hand, since $V$ is $L$-stable, $L\cdot (2G_0 - D)< 0$. Since $L$ and $M$ are
not separated by any wall of type $(w,p)$, it follows that $M\cdot (2G_0 -
D)\leq 0$. Finally using
$$\align
p_1(\operatorname{ad}W) &= (2G_0 - D +C)^2 - 4\ell (\Cal Z_0) \\
&= p_1(\operatorname{ad}V) + 2D\cdot C + (C)^2 -
4f\\
&> p + k^2(M)^2,
\endalign$$
we obtain
$$(2G_0 - D)^2 +2k(2G_0 - D)\cdot M > p.$$
Let $m = -(2G_0 - D)\cdot M$. As we have seen above $m\leq kM^2$ and $m\geq 0$.
The above inequality can be rewritten as
$$ 2km < (2G_0 - D)^2 -p.$$
We claim that $m=0$. Otherwise
$$2k < \frac{(2G_0 - D)^2}{m}  -\frac{p}{m} .$$
By the Hodge index theorem $(2G_0 - D)^2 M^2 \leq \left[(2G_0 - D)\cdot
M\right]^2 = m^2$, so that $(2G_0 - D)^2 \leq m^2/M^2$. Plugging this into the
inequality above, using $-p\geq 0$, gives
$$2k < \frac{m}{M^2} - \frac{p}{m} \leq k - p,$$
i\.e\. $k< -p$, contradicting our choice of $k$. Thus $m= -(2G_0 - D)\cdot M
=0$.

Now the inclusions $\scrO_X(G_0) \subset W\subset V$ define an inclusion
$\scrO_X(G_0) \subset V$. Thus there is an effective divisor $E$ and an
inclusion $\scrO_X(G_0+E) \to V$ with torsion free cokernel. Let $G = G_0 +E$.
Thus there is an exact sequence
$$0 \to \scrO_X(G) \to V \to \scrO_X(-G+D)\otimes I_{\Cal Z}\to 0.$$
We claim that $(2G-D)\cdot M = 0$. Since $V$ is $L$-stable, $(2G-D)\cdot L <
0$, and since $L$ and $M$ are not separated by a wall of type $(w,p)$,
$(2G-D)\cdot M \leq 0$. On the other hand,
$$(2G-D)\cdot M = (2G_0 - D) \cdot M
+ 2(E\cdot M) = -m + 2(E\cdot M)=2(E\cdot M).$$
As $E$ is effective and $M$ is nef, $2(E\cdot M) \geq 0$. Thus $(2G-D)\cdot M =
0$. As $M^2 >0$, we must have $(2G-D)^2<0$. Using $p = (2G-D)^2 -4\ell(\Cal
Z)\leq (2G-D)^2$, we see that $2G-D$ is a wall of type $(w,p)$.

Finally note that $\operatorname{Supp}\Cal Z\cap C \neq \emptyset$, for
otherwise we would have $V|C$ semistable. This concludes the proof of Lemma
2.2.
\endproof

Returning to the proof of Lemma 2.1, we claim first that, given $k\gg 0$ and
$C\in |kM|$, for all $N$ sufficiently large the sections of
$\Cal L (C, \theta)^N$ define a base point free linear series on
$\overline{\frak M_L(w,p)}$. To see this, we first claim that, for $k \gg 0$,
and for a generic $C\in |kM|$, the restriction map $V\mapsto V|C$ defines a
rational map $r_C\: \frak M_L(w,p)\dasharrow \frak M(C)$, where $\frak M(C)$ is
the moduli space of equivalence classes of semistable rank two bundles on $C$
such that the parity of the determinant is even. It suffices to prove that, for
every component $N$ of $\frak M_L(w,p)$ there is one $V\in N$ and
one $C \in |kM|$ such that $V|C$ is semistable, for then the same will hold for
a Zariski open subset of $|kM|$. Now given $V$, choose a fixed $C_0 \in |kM|$.
If $V|C_0$ is not semistable, then by Lemma 2.2 there is an exact sequence
$$0 \to \scrO_X(G) \to V \to \scrO_X(-G+D)\otimes I_{\Cal Z} \to 0,$$
where $\Cal Z$ is a zero-dimensional subscheme of $X$ meeting $C_0$. Choosing
$C$ to be a curve in $|kM|$ disjoint from $\Cal Z$, which is possible since
$|kM|$ is base point free, it follows that the restriction $V|C$ is semistable.

For $C$ fixed, let
$$\align
B_C= \{\,V \in \overline{\frak M_L(w,p)}:&\text{ either $V$ is not locally free
over some point of $C$}\\
&\text{ or  $V|C$  is not semistable }\,\}.
\endalign$$
By the openness of stability and local freeness, the set $B_C$ is a closed
subset of $\overline{\frak M_L(w,p)}$ and
$r_C$ defines a morphism from
$\overline{\frak M_L(w,p)} -B_C$ to
$\frak M(C)$. Standard estimates (cf\. [10], [12], [32], [24], [27]) show
that, possibly after modifying the constant $A$ introduced at the beginning of
the proof of Theorem 1.4, the codimension of
$B_C$ is at least two in $\overline{\frak M_L(w,p)}$ provided that $p\leq A$
(where as usual $A$ is independent of $k$ and depends only on $X$ and $M$).
Indeed the set of bundles
$V$ which fit into an exact sequence
$$0 \to \scrO_X(G) \to V \to \scrO_X(D-G)\otimes I_{\Cal Z} \to 0,$$
where $G$ is a divisor such that $(2G-D)\cdot M = 0$,
may be parametrized by a scheme of dimension $-\frac34p + O(\sqrt{|p|})$ by
e\.g\. [12], Theorem 8.18. Moreover the constant implicit in the notation
$O(\sqrt{|p|})$ can be chosen uniformly over $\Cal N$. The case of
nonlocally free
$V$ is taken care of by assumption (2) in the discussion of the constant $A$:
it follows from standard deformation theory (see again [12], [24]) that at a
generic point of the locus of nonlocally free sheaves corresponding to the
semistable torsion free sheaf $V$ the deformations of $V$ are versal for the
local  deformations of the  singularities of $V$. Thus for a general nonlocally
free $V$,
$V$ has just one singular point which is at a general point of $X$ and so does
not lie on $C$. Thus the set of
$V$ which are not locally free at some point of $C$ has codimension at least
two
(in fact exactly two) in $\overline{\frak M_L(w,p)}$.

Let
$\Cal L_C$ be the determinant line bundle on $\frak M(C)$ associated to the
line bundle $\theta$ (see for instance [15] Chapter 5 Section 2). Then by
definition the pullback via
$r_C$ of
$\Cal L_C$ is the restriction of $\Cal L (C, \theta)$ to $\overline{\frak
M_L(w,p)} -B_C$. Since $B_C$ has codimension two, the sections of
$\Cal L_C^N$ pull back to sections of $\Cal L (C, \theta)^N$ on
$\overline{\frak M_L(w,p)}$. Since $\Cal L_C$ is ample, given $V\in
\overline{\frak M_L(w,p)} -B_C$, there exists an $N$ and a section of $\Cal
L_C^N$ not vanishing at
$r_C(V)$, and thus there is a section of $\Cal L (C, \theta)^N$ not vanishing
at $V$. Moreover by [23], for all smooth $C'\in |kM|$ and choice of an
appropriate line bundle $\theta '$ on $C'$, there is an isomorphism $\Cal L (C,
\theta)^N \cong
\Cal L(C', \theta ')^N$. Next we claim that, for every $V\in
\overline{\frak M_L(w,p)}$, there exists a $C$ such that $V$ is locally free
above $C$ and $V|C$ is semistable. Given $V$, it fails to be locally free
at a finite set of points, and its double dual $W$ is again semistable. Thus
applying the above to $W$, and again using the fact that $|kM|$ has no base
points, we can find $C$ such that $V$ is locally free over $C$ and such that
$V|C = W|C$ is semistable. Thus, given $V$, there exists an $N$ and a section
of $\Cal L (C, \theta)^N$ which does not vanish at $V$. Since $\overline{\frak
M_L(w,p)}$ is of finite type, there exists an $N$ which works for all $V$, so
that the linear system corresponding to $\Cal L (C, \theta)^N$ has no base
points.

Finally we must show that, for $k\gg 0$, the morphism induced by $\Cal L (C,
\theta)^N$ is in fact generically finite for $N$ large. We claim that it
suffices to show that the restriction of the rational map $r_C$ to
$\overline{\frak M_L(w,p)} -B_C$ is generically finite (it is here that we must
use the condition  on the singularities of
$\bar X$ in the statement of Theorem 1.4). Supposing this to be the case, and
fixing a $V \in \overline{\frak M_L(w,p)} -B_C$ for which
$r_C^{-1}(r_C(V))$ is finite, we consider the intersection of all the divisors
in $\Cal L (C, \theta)^N$ containing $V$, where $N$ is chosen so that $\Cal
L_C^N$ is very ample. This intersection always contains $V$ and is a subset of
$r_C^{-1}(r_C(V)) \cup B_C$. In particular $V$ is an isolated point of the
fiber, and so the morphism defined by $\Cal L (C, \theta)^N$ cannot have all
fibers of purely positive dimension. Thus it is generically finite.

To see that $r_C$ is generically finite, we shall show that, for generic
$V$, the restriction map
$$r\: H^1(X;\operatorname{ad}V) \to H^1(C;\operatorname{ad}V|C))$$
is injective. The map $r$ is just the differential of the map $r_C$ from $\frak
M_L(w,p)$ to $\frak M(C)$ at the point corresponding to $V$, and so if $V$ is
generic then $r_C$ is finite. Now the kernel of the map
$r$ is a quotient of
$H^1(X;  \operatorname{ad}V\otimes \scrO_X(-C))$, and we need to find
circumstances where this group is zero, at least if
$C\in |kM|$ for
$k$ sufficiently large. By Serre duality it suffices to show that $H^1(X;
\operatorname{ad}V\otimes \scrO_X(C)\otimes K_X)=0$ for $k$ sufficiently large.
By applying the Leray spectral sequence to the morphism
$\varphi\: X\to \bar X$, it suffices to show that
$$H^1(\bar X;R^0\varphi _* (\operatorname{ad}V\otimes \scrO_X(C)\otimes
K_X)) =0$$
and that
$R^1\varphi _* (\operatorname{ad}V\otimes \scrO_X(C)\otimes
K_X)=0$. Now $M$ is the pullback of an ample line bundle $\bar M$ on
$\bar X$, and $\scrO_X(C)$ is the pullback of $(\bar M)^{\otimes k}$. Thus for
fixed $V$ and $k \gg 0$,
$$\align
&H^1(\bar X;R^0\varphi _* (\operatorname{ad}V\otimes \scrO_X(C)\otimes
K_X)) \\=
&H^1(\bar X;R^0\varphi _* (\operatorname{ad}V\otimes
K_X)\otimes (\bar M^k)) =0.
\endalign$$
Moreover $R^1\varphi _* (\operatorname{ad}V\otimes \scrO_X(C)\otimes
K_X)) = R^1\varphi _* (\operatorname{ad}V\otimes
K_X)\otimes (\bar M^k)$, so that it is enough to show that $R^1\varphi _*
(\operatorname{ad}V\otimes K_X)=0$. By the formal functions theorem,
$$R^1\varphi _* (\operatorname{ad}V\otimes K_X) = \varprojlim _mH^1(mZ;
\operatorname{ad}V\otimes K_X |mZ),$$ where $Z = \bigcup Z_i$ is the union of
the connected components
$Z_i$ of the one-dimensional fibers of
$\varphi$. Thus it suffices to show that, for all $i$ and all positive integers
$m$,
$H^1(mZ_i; \operatorname{ad}V\otimes K_X |mZ_i) = 0$. Now by the adjunction
formula
$\omega _{mZ_i} = K_X\otimes \scrO_X(mZ_i)|mZ_i$, where $\omega _{mZ_i}$ is the
dualizing sheaf of the Gorenstein scheme $mZ_i$. Thus $K_X|mZ_i =
\scrO_X(-mZ_i) |mZ_i\otimes\omega _{mZ_i}$ and we must show the vanishing of
$$H^1(mZ_i; (\operatorname{ad}V\otimes \scrO_X(-mZ_i)) |mZ_i\otimes\omega
_{mZ_i}).$$ By Serre duality, it suffices to show that, for all $m>0$,
$$H^0(mZ_i; (\operatorname{ad}V\otimes \scrO_X(mZ_i)) |mZ_i)=0.$$
We shall deal with this problem in the next subsection.
\medskip

\noindent {\bf Remark.} (1) Instead of arguing that the restriction map $r_C$
was generically finite, one could also check that it was generically one-to-one
by showing that for generic $V_1$, $V_2$, the restriction map
$$H^0(X; Hom (V_1, V_2)) \to H^0(C; Hom (V_1, V_2)|C)$$
is surjective (since then an isomorphism from $V_1|C$ to $V_2|C$ lifts to a
nonzero map from $V_1$ to $V_2$, necessarily an isomorphism by stability). In
turn this would have amounted to showing that
$H^1(X;  Hom(V_1, V_2)\otimes \scrO_X(-C))=0$ for generic $V_1$ and $V_2$, and
this would have been essentially the same argument.

\smallskip
(2) Suppose that $\varphi \: X\to \bar X$ is the blowup of a smooth surface
$\bar X$ at a point $x$, and that $M$ is the pullback of an ample divisor on
$\bar X$. Let
$Z\cong\Pee ^1$ be the exceptional curve. In this case, if $c_1(V)\cdot Z$ is
odd, say $2a+1$, then the generic behavior for $V|Z$ is $V|Z \cong \scrO_{\Pee
^1}(a)\oplus
\scrO_{\Pee ^1}(a+1)$ and the restriction map
exhibits $\frak M_L(w,p)$ (generically) as a $\Pee ^1$-bundle over its image
(see for instance [5]). Thus the hypothesis that $\varphi$ contracts no
exceptional curve is essential.

\ssection{2.3.  Restriction of stable bundles to certain curves}

Let us recall the
basic properties of rational and minimally elliptic singularities. Let $x$ be a
normal singular point on a complex surface $\bar X$, and let $\varphi \: X \to
\bar X$ be the minimal resolution of singularities of $\bar X$. Supppose that
$\varphi ^{-1}(x) = \bigcup _iD_i$. The singularity is a {\sl rational\/}
singularity if
$(R^1\varphi _*\scrO_X)_x = 0$. Equivalently, by [1], $x$ is rational if and
only if, for every choice of nonnegative integers $n_i$ such that at least one
of the $n_i$ is strictly positive, if we set $Z = \sum _in_iD_i$, the
arithmetic genus $p_a(Z)$ of the effective curve $Z$ satisfies $p_a(Z) \leq 0$.
Here $p_a(Z) = 1 - \chi (\scrO_Z) = 1 - h^0(\scrO_Z) + h^1(\scrO_Z)\leq
h^1(\scrO_Z)$; moreover we have the adjunction formula
$$p_a(Z) = 1 + \frac12 (K_X+Z)\cdot Z.$$
Now every minimal resolution of a normal surface singularity $x$ has a {\sl
fundamental cycle\/} $Z_0$, which is an effective cycle $Z_0$ supported in the
set $\varphi ^{-1}(x)$ and satisfying $Z_0
\cdot D_i \leq 0$ and $Z_0 \cdot D_i < 0$ for some $i$ which is minimal with
respect to the above properties. We may find $Z_0$ as follows [22]: start
with an arbitrary component $A_1$ of $\varphi ^{-1}(x)$ and set $Z_1= A_1$.
Now either $Z_0 = A_1$ or there exists another component $A_2$ with $Z_1\cdot
A_2>0$. Set $Z_2 = Z_1 + A_2$ and continue this process. Eventually we reach
$Z_k = Z_0$. Such a sequence $A_1, \dots , A_k$ with $Z_i = \sum _{j\leq
i}A_j$ and $Z_i\cdot A_{i+1} >0$, $Z_k = Z_0$ is called a {\sl computation
sequence}. By a theorem of Artin [1],
$x$ is rational if and only if $p_a(Z_0) \leq 0$, where $Z_0$ is the
fundamental
cycle, if and only if $p_a(Z_0) = 0$. Moreover, if
$x$ is a rational singularity, then every component $D_i$ of $\varphi ^{-1}(x)$
is a smooth rational curve, the
$D_i$ meet transversally at at most one point, and the dual graph of $\varphi
^{-1}(x)$ is contractible.

Next we recall the properties of minimally elliptic singularities [22]. A
singularity $x$ is {\sl minimally elliptic\/} if and only if there exists a
{\sl minimally elliptic cycle\/} $Z$ for $x$, in other words a cycle $Z= \sum
_in_iD_i$ with all $n_i>0$ such that $p_a(Z) = 1$ and $p_a(Z') \leq 0$ for all
nonzero effective cycles $Z'<Z$ (i\.e\. such that $Z' = \sum _in_i'D_i$ with
$0\leq n_i'\leq n_i$ and $Z'\neq Z$). In this case it follows that $Z=Z_0$ is
the fundamental cycle for $x$, and $(K_X+Z_0)\cdot D_i = 0$ for every component
$D_i$ of
$\varphi ^{-1}(x)$. If $Z_0$ is reduced, i\.e\. if $n_i = 1$ for all $i$, then
the possibilities for $x$ are as follows:
\roster
\item $\varphi ^{-1}(x)$ is an irreducible curve of arithmetic genus one, and
thus is either a smooth elliptic curve or a singular rational curve with
either a node or a cusp;
\item $\varphi ^{-1}(x)= \bigcup _{i=1}^tD_i$ is a cycle of $t\geq 2$ smooth
rational curves meeting transversally, i\.e\.  $D_i\cdot D_{i+1} =1$, $D_i
\cdot D_j \neq 0$ if and only if $i \equiv j \pm 1 \mod t$, except for $t=2$
where $D_1 \cdot D_2 = 2$;
\item $\varphi ^{-1}(x)= D_1\cup D_2$, where the $D_i$ are smooth
rational, $D_1\cdot D_2 = 2$ and
$D_1\cap D_2 $ is a single point (so that $\varphi ^{-1}(x)$ has a tacnode
singularity) or $\varphi ^{-1}(x)= D_1\cup D_2 \cup D_3$ where the $D_i$ are
smooth rational, $D_i\cdot D_j =
1$ but
$D_1\cap D_2 \cap D_3$ is a single point (the three curves meet at a common
point).
\endroster
Here $x$ is called a {\sl simple elliptic singularity\/} in case $\varphi
^{-1}(x)$ is a smooth elliptic curve, a {\sl cusp singularity\/} if $\varphi
^{-1}(x)$ is an irreducible rational curve with a node or a cycle as in (2),
and a {\sl triangle singularity\/} in the remaining cases.
If $Z_0$ is not reduced, then all components $D_i$ of $\varphi ^{-1}(x)$ are
smooth rational curves meeting transversally and the dual graph of $\varphi
^{-1}(x)$ is contractible.

With this said, and using the discussion in the previous subsection, we will
complete the proof of Theorem 1.4 by showing that
$H^0(mZ_i; \operatorname{ad}V\otimes \scrO_X(mZ_i)|mZ_i) = 0$ for all $i$,
where $x_1, \dots, x_k$ are the singular points of $\bar X$ and $Z_i$ is an
effective cycle with $\operatorname{Supp}Z_i = \varphi ^{-1}(x_i)$. The precise
statement is as follows:

\theorem{2.3} Let $\varphi \: X \to \bar X$ be a birational morphism from $X$
to a normal projective surface $\bar X$, corresponding to a nef, big, and
eventually base point free divisor  $M$. Let $w$ be the \rom{mod} $2$ reduction
of $[K_X]$, and suppose that
\roster
\item"{(i)}" $\varphi$ contracts no exceptional curve; in other words, if $E$
is an exceptional curve of the first kind on $X$, then $M\cdot E>0$.
\item"{(ii)}" $\bar X$ has only rational and minimally elliptic singularities.
\endroster
Then there exists a  constant $A$ depending only
on $p$ and $\Cal N$ with the following property: for every singular point
$x$ of $\bar X$, there exists an effective cycle $Z$ with
$\operatorname{Supp}Z = \varphi^{-1}(x)$  such that, for all ample line
bundles $L$ in
$\Cal N$, all
$p$ with $p \leq A$, and generic bundles $V$ in $\frak M_L(w,p)$,
$$H^0(mZ; \operatorname{ad}V\otimes \scrO_X(mZ)|mZ) = 0$$
for every positive integer $m$.
\endstatement

The statement of (i) may be rephrased by saying that $X$ is the {\sl minimal
resolution\/} of $\bar X$. As $\operatorname{ad}V \subset Hom(V,V)$,
it suffices to prove that $H^0(mZ; Hom(V,V)\otimes \scrO_X(mZ)|mZ) = 0$. We
will consider the case of rational singularities and minimally elliptic
singularities separately. Let us begin with the proof for rational
singularities. Let $\varphi ^{-1}(x)  = \bigcup _iD_i$, where each
$D_i$ is a smooth rational curve. By the assumption (4) of the previous
subsection, we can assume that the constant $A$ has been chosen so that
$V|D_i$ is a generic bundle over $D_i\cong \Pee ^1$ for every
$i$. Thus either there exists an integer $a$ such that $V|D_i
\cong \scrO_{\Pee^1}(a)
\oplus \scrO_{\Pee^1}(a+1)$, if $w\cdot D_i \neq 0$, or there exists an $a$
such that $V|D_i  \cong \scrO_{\Pee^1}(a) \oplus
\scrO_{\Pee^1}(a)$, if $w\cdot D_i = 0$. Next, we have the following claim:

\claim{2.4} Suppose that $x$ is a rational singularity. Let $\varphi\:X \to
\bar X$ be a resolution of $x$. There exist a sequence of curves $B_0, \dots,
B_k$, such that $B_i \subseteq \varphi ^{-1}(x)$ for all $i$, with the
following property:
\roster
\item Let $C_i = \sum _{j\leq i}B_i$. Then $B_i \cdot C_i \leq B_i^2 + 1$.
\item $C_k = Z_0$, the fundamental cycle of $x$.
\endroster
\endstatement
\proof Since $(K_X+Z_0)\cdot Z_0<0$, there must exist
a component $B^{(0)} = D_i$ of $\operatorname{Supp}Z_0 = \varphi ^{-1}(x)$ such
that $(K_X+Z_0)\cdot B^{(0)} < 0$. Thus
$$Z_0 \cdot B^{(0)} < -K_X\cdot B^{(0)} = (B^{(0)})^2 + 2.$$
Set $Z_1 = Z_0 - B^{(0)}$. Suppose that $Z_1$ is nonzero. Then $Z_1$ is again
effective, and by Artin's criterion $p_a(Z_1) \leq 0$. Thus by repeating the
above argument there is a
$B^{(1)}$ contained in the support of $Z_1$ such that $Z_1 \cdot B^{(1)} <
(B^{(1)})^2 + 2$. Continuing, we eventually find $B^{(2)}, \dots, B^{(k)}$
with $B^{(i)}$ contained in the support of $Z_i$, $Z_{i+1} = Z_i -B^{(i)}$ and
$Z_k = B^{(k)}$, and such that $Z_i \cdot B^{(i)} < (B^{(i)})^2 +2$. If we now
relabel $B^{(i)} = B_{k-i}$, then $Z_i = \sum _{j\leq n-i}B_j$ and the curves
$B_0,
\dots, B_k$ are as claimed.
\endproof

Returning to the proof of Theorem 2.3, we first prove that
$$H^0(Z_0; Hom(V, V)\otimes \scrO_X(Z_0)|Z_0) = 0.$$
We have the exact sequence
$$0 \to \scrO_{C_{i-1}}(C_{i-1}) \to \scrO_{C_i}(C_i) \to \scrO_{B_i}(C_i) \to
0.$$
Tensor this sequence by $Hom (V, V)$. We shall prove by induction that
$$H^0(Hom (V, V) \otimes \scrO_{C_i}(C_i)) = 0$$ for all $i$. It suffices to
show that
$H^0(Hom (V, V) \otimes \scrO_{B_i}(C_i)) = 0$ for all $i$. Now
$\scrO_{B_i}(C_i)$ is a line bundle on the smooth rational curve $B_i$. If
$V|B_i \cong \scrO_{\Pee^1}(a)
\oplus \scrO_{\Pee^1}(a+1)$, then $w\cdot B_i \neq 0$ and so $B_i^2$ is
odd. Since $B_i$ is not an exceptional curve, $B_i^2
\leq -3$ and so
$B_i \cdot C_i \leq -2$. Thus, as
$$Hom (V, V)|B_i = \scrO_{\Pee^1}(-1) \oplus \scrO_{\Pee^1}
\oplus \scrO_{\Pee^1} \oplus \scrO_{\Pee^1}(1),$$
we see that $H^0(Hom (V, V)\otimes \scrO_{B_i}(C_i)) = 0$. Likewise if
$V|D_i  \cong \scrO_{\Pee^1}(a) \oplus \scrO_{\Pee^1}(a)$, then
using
$B_i \cdot C_i \leq -1$ we again have $H^0(Hom (V, V) \otimes
\scrO_{B_i}(C_i)) = 0$. Thus by induction
$$H^0(Hom (V, V) \otimes \scrO_{C_k}(C_k)) = H^0(Hom (V, V) \otimes
\scrO_{Z_0}(Z_0))=0.$$
The vanishing of $H^0(mZ_0; Hom(V, V)\otimes \scrO_X(mZ_0)|mZ_0)$ is
similar, using instead the exact sequence
$$0 \to \scrO_{mZ_0+C_{i-1}}(mZ_0 + C_{i-1}) \to \scrO_{mZ_0 +C_i}(mZ_0 + C_i)
\to \scrO_{B_i}(mZ_0 + C_i) \to 0.$$
This concludes the proof in the case of a rational singularity.

For minimally elliptic
singularities, we shall deduce the theorem from the following more general
result:

\theorem{2.5} Let $\varphi \: X \to \bar X$ be a birational morphism from $X$
to a normal projective surface $\bar X$, corresponding to a nef, big, and
eventually base point free divisor  $M$. Let $w$ be an arbitrary element of
$H^2(X; \Zee/2\Zee)$, and suppose that
\roster
\item"{(i)}" $\varphi$ contracts no exceptional curve; in other words, if $E$
is an exceptional curve of the first kind on $X$, then $M\cdot E>0$.
\item"{(ii)}" If $D$ is a component of $\varphi^{-1}(x)$ such that $w\cdot D
\neq 0$, then $Z_0 \cdot D <0$,
where $Z_0$ is the fundamental cycle of $\varphi^{-1}(x)$.
\endroster
Then the conclusions of Theorem \rom{2.3} hold for the moduli space $\frak
M_L(w,p)$ for all $p\ll 0$. In particular the conclusions of Theorem \rom{2.3}
hold if
$\varphi^{-1}(x)$ is irreducible.
\endstatement

\demo{Proof that \rom{(2.5)} implies \rom{(2.3)}} We must show that every
minimally elliptic singularity satisfies the hypotheses of Theorem 2.5(ii),
provided that
$w$ is the mod 2 reduction of $K_X$. Suppose that $x$ is minimally elliptic and
that $w\cdot D \neq 0$. Thus $K_X\cdot D$ is odd. Moreover if $D$ is smooth
rational then $D^2\neq -1$ and $K_X\cdot D \geq 0$ so that $K_X\cdot D \geq 1$.
Now $(K_X+Z_0)\cdot D = 0$. Thus $Z_0
\cdot D = -(K_X\cdot D)\leq -1$. Likewise if $p_a(D) \neq 0$, so that $D$ is
not a smooth rational curve, then
$\varphi ^{-1}(x) = D$ is an irreducible curve and (2.3) again follows.
\endproof

\demo{Proof of Theorem \rom{2.5}} We begin with a  lemma on sections
of line bundles over effective cycles supported in $\varphi^{-1}(x)$, which
generalizes (2.6) of [22]:

\lemma{2.6}
Let $Z_0$ be the fundamental cycle of $\varphi ^{-1}(x)$ and let $\lambda$
be a line bundle on $Z_0$ such that $\deg (\lambda|D)\leq 0$ for each component
$D$ of the support of $Z_0$. Then either $H^0(Z_0; \lambda) =0$ or $\lambda =
\scrO_{Z_0}$ and $H^0(Z_0; \lambda) \cong \Cee$.
\endstatement
\proof Choose a computation sequence for $Z_0$, say $A_1, A_2, \dots, A_k$.
Thus, if we set $Z_i = \sum _{j\leq i}A_j$, then $Z_i\cdot A_{i+1} >0$, and
$Z_k = Z_0$. Now we have an exact sequence
$$0 \to \scrO_{A_{i+1}}(-Z_i)\to \scrO_{Z_{i+1}} \to \scrO_{Z_i} \to 0.$$
Thus $\deg (\scrO_{A_{i+1}}(-Z_i) \otimes \lambda|A_{i+1})<0$. It follows that
$H^0(\scrO_{Z_{i+1}}\otimes \lambda) \subseteq H^0(\scrO_{Z_i}\otimes
\lambda)$ for all $i$. By induction $\dim H^0(\scrO_{Z_i}\otimes
\lambda) \leq 1$ for all $i$, $1\leq i \leq k$. Thus  $\dim H^0(Z_0; \lambda)
\leq 1$. Moreover, if $\dim H^0(Z_0; \lambda) = 1$, then the natural map
$$H^0(\scrO_{Z_{i+1}}\otimes \lambda) \to H^0(\scrO_{Z_i}\otimes
\lambda)$$ is an isomorphism for all $i$, and so the induced map $H^0(Z_0;
\lambda) \to H^0(A_1; \lambda |A_1)$ is an isomorphism and $\dim H^0(A_1;
\lambda |A_1) =1$. Thus $\lambda |A_1$ is trivial and a nonzero section of
$H^0(Z_0; \lambda)$ restricts to a generator of $\lambda |A_1$. Since we can
begin a computation sequence with an arbitrary choice of $A_1$, we see that a
nonzero section $s$ of $H^0(Z_0; \lambda)$ restricts to a nonvanishing section
of $H^0(D; \lambda |D)$ for every $D$ in the support of $\varphi ^{-1}(x)$.
Thus
the map $\scrO_{Z_0} \to \lambda$ defined by $s$ is an isomorphism.
\endproof

\noindent {\bf Remark.} The lemma is also true if $\lambda$ is allowed to have
degree one on some components $D$ of $Z_0$ with $p_a(D)\geq 2$, provided that
$\lambda |D$ is general for these components, and a slight variation holds if
$\lambda$ is also allowed to have degree one on some components $D$ of $Z_0$
with $p_a(D)=1$.
\medskip

We next construct a bundle $W$ over $Z_0$ with certain vanishing properties:

\lemma{2.7} Suppose that $\varphi \: X \to \bar X$ is the minimal resolution of
the normal surface singularity $x$. Let
$\mu$ be a line bundle over the scheme
$Z_0$. Suppose further that, if $D$ is a component of $\varphi^{-1}(x)$ such
that $\deg (\mu |D)$ is odd, then
$Z_0 \cdot D <0$, where $Z_0$ is the fundamental cycle of $\varphi^{-1}(x)$.
Then there exists a rank two vector bundle
$W$ over $Z_0$ with $\det W = \mu$ and such that
$$H^0(Z_0; Hom(W,W)\otimes \scrO_X(mZ_0)|Z_0) = 0$$
for every $m\geq 1$.
\endstatement
\proof Let $\varphi ^{-1}(x) = \bigcup_iD_i$. Then there exists an
integer $a_i$ such that $\deg\mu |D_i = 2a_i$ or
$2a_i+1$, depending on whether $\deg (\mu| D_i)$ is odd or even. Since $\dim
Z_0 = 1$, the natural maps $\operatorname{Pic}Z_0 \to
\operatorname{Pic}(Z_0)_{\text{red}} \to \bigoplus _i\operatorname{Pic}D_i$ are
surjective. Thus we may choose a line bundle $L_1$ over $Z_0$ such that $\deg
(L_1|D_i) = a_i$. It follows that $\mu \otimes L_1^{\otimes
- -2}|D_i$ is a line bundle over $D_i$ of degree zero or 1, and if it is of
degree 1, then $Z_0\cdot D_i <0$. Hence $\mu \otimes
L_1^{\otimes -2}\otimes \scrO_{Z_0}(Z_0)$ has degree at most zero on $D_i$ for
every $i$.

Set $L_2 = \mu \otimes L_1^{-1}$. Thus $L_1 \otimes L_2 = \mu$ and $\deg
(L_2|D_i) = a_i$ or $a_i +1$ depending on whether $\deg (\mu| D_i)$ is even or
odd. The line bundle $L_1^{-1}\otimes L_2=\mu \otimes
L_1^{\otimes -2}$ thus has degree zero on
those components $D_i$ such that $\deg (\mu| D_i)$ is even and 1 on the
components $D_i$ such that $\deg (\mu| D_i)$ is odd. Moreover $\deg
(L_1^{-1}\otimes L_2 \otimes \scrO_{Z_0}(mZ_0)|D_i) \leq 0$ for every $i$.

\claim{2.8} Under the assumptions of \rom{(2.7)}, there exists a
nonsplit extension
$W$ of $L_2$ by $L_1$ except in the case where $x$ is rational, $\deg
(\mu| D_i)$ is odd for at most one $i$, and the multiplicity of $D_i$ in $Z_0$
is one for such $i$, or $\chi (\scrO_{Z_0}) = 0$ and $\deg \mu |D_i$ is even
for
every $i$.
\endstatement
\proof  A nonsplit extension exists if and only if  $h^1(L_2^{-1}\otimes L_1)
\neq 0$. Now by the Riemann-Roch theorem applied to
$Z_0 =
\sum _in_iD_i$, we have
$$h^1(Z_0; L_2^{-1}\otimes L_1) = h^0(Z_0; L_2^{-1}\otimes L_1) -\sum _in_i\deg
(L_2^{-1}\otimes L_1|D_i) - \chi (\scrO_{Z_0}).$$
Here $\deg (L_2^{-1}\otimes L_1|D_i)= 0$ on those $D_i$ with $\deg (\mu| D_i)$
even and $=-1$ on the $D_i$ with $\deg (\mu| D_i)$ odd.  Moreover
$h^0(\scrO_{Z_0}) = 1$ by Lemma 2.6 and so  $\chi (\scrO_{Z_0})\leq 1$, with
$\chi (\scrO_{Z_0}) =1$ if and only if $x$ is rational. Thus
$$h^1(Z_0; L_2^{-1}\otimes L_1) \geq \sum \{\, n_i\: \deg (\mu| D_i) \text{ is
odd}\,\} - \chi (\scrO_{Z_0}).$$
Hence if $h^1(Z_0; L_2^{-1}\otimes L_1) =0$, then either $x$ is rational, $\deg
(\mu| D_i)$ is odd for at most one $i$, and for such $i$ the multiplicity of
$D_i$ in $Z_0$ is one, or $\deg (\mu |D_i)$ is even for all $i$ and $\chi
(\scrO_{Z_0}) = 0$.
\endproof

Returning to the proof of (2.7), choose $W$ to be a nonsplit extension of $L_2$
by $L_1$ if such exist, and set $W = L_1 \oplus L_2$ otherwise. To see that
$H^0(Z_0; Hom(W,W)\otimes \scrO_X(mZ_0)|Z_0) = 0$, we consider the two exact
sequences
$$\gather
0 \to L_1 \to W \to L_2 \to 0;\\
0 \to L_1 \otimes \scrO_{Z_0}(mZ_0)\to W\otimes \scrO_{Z_0}(mZ_0) \to
L_2\otimes \scrO_{Z_0}(mZ_0) \to 0.
\endgather$$
Clearly $H^0(Z_0; Hom(W,W)\otimes \scrO_X(mZ_0)|Z_0) = 0$ if
$$H^0(L_1^{-1}\otimes L_2 \otimes \scrO_{Z_0}(mZ_0)) = H^0(\scrO_{Z_0}(mZ_0)) =
H^0(L_2^{-1}\otimes L_1 \otimes \scrO_{Z_0}(mZ_0)) =0.$$
The line bundles $\scrO_{Z_0}(mZ_0)$ and $L_2^{-1}\otimes L_1 \otimes
\scrO_{Z_0}(mZ_0)$ have nonpositive degree on each $D_i$ and (since $Z_0 \cdot
D_i <0$ for some
$i$) have strictly negative degree on at least one component. Thus by Lemma 2.6
$H^0(\scrO_{Z_0}(mZ_0))$ and $H^0(L_2^{-1}\otimes L_1 \otimes
\scrO_{Z_0}(mZ_0))$ are both zero. Let us now consider the group
$H^0(L_1^{-1}\otimes L_2 \otimes \scrO_{Z_0}(mZ_0))$. By the hypothesis that
$Z_0\cdot D_i < 0$ for each $D_i$ such that $\deg (\mu |D_i)$ is odd, the line
bundle $L_1^{-1}\otimes L_2 \otimes \scrO_{Z_0}(mZ_0)$ has nonpositive degree
on all components $D_i$. Thus by Lemma 2.6 either $H^0(L_1^{-1}\otimes L_2
\otimes \scrO_{Z_0}(mZ_0)) = 0$ or $L_1^{-1}\otimes L_2 \otimes \scrO_{Z_0}
(mZ_0) \cong \scrO_{Z_0}$. Clearly this last case is only possible if $m=1$ and
$L_1 \cong L_2 \otimes \scrO_{Z_0}(Z_0)$, and if moreover $Z_0 \cdot D_i =0$ if
$\deg(\mu |D_i)$ is even and $Z_0 \cdot D_i =-1$ if $\deg(\mu
|D_i)$ is odd. As $Z_0\cdot D_i <0$ for at least one $i$, $\deg(\mu
|D_i)$ is odd for at least one $i$ as well. In this case, if the nonzero
section of $L_1^{-1}\otimes L_2 \otimes \scrO_{Z_0}(Z_0)$ lifts to give a map
$L_1\to W\otimes \scrO_{Z_0}(Z_0)$, then the image of
$L_1$ in $W\otimes \scrO_{Z_0}(Z_0)$ splits the exact sequence
$$0 \to L_1 \otimes \scrO_{Z_0}(Z_0)\to W\otimes \scrO_{Z_0}(Z_0) \to
L_2\otimes \scrO_{Z_0}(Z_0) \to 0.$$
Thus $W$ is also a split extension. By Claim 2.8, since $\deg (\mu |D_i)$ is
odd for at least one $i$, it must therefore be the case that $x$ is rational,
$\deg (\mu| D_i)$ is odd for exactly one $i$, and for such
$i$ the multiplicity of $D_i$ in $Z_0$ is one. Moreover $Z_0\cdot D_j \neq 0$
exactly when $j=i$ and in this case $Z_0 \cdot D_i = -1$. But as the
multiplicity of
$D_i$ in $Z_0$ is 1, we can write $Z_0 = D_i + \sum _{j\neq i}n_jD_j$, and thus
$$Z_0^2 = Z_0\cdot D_i = -1.$$
By a theorem of Artin [1], however, $-Z_0^2$ is the multiplicity of the
rational
singularity $x$. It follows that $x$ is a smooth point and $\varphi$ is the
contraction of a generalized exceptional curve, contrary to hypothesis. This
concludes the proof of (2.7).
\endproof

We may now finish the proof of (2.5). Start with a generic vector bundle $V_0
\in \frak M_L(w,p)$ on
$X$ satisfying the condition that $H^1(X; \operatorname{ad}V_0)\to H^1(Z_0;
\operatorname{ad}V_0|Z_0)$ is surjective. If $\mu = \det V_0|Z_0$, note that,
according to the assumptions of (2.5),
$\mu$ satisfies the hypotheses of Lemma 2.7. For $V\in \frak M_L(w,p)$, let
$$H(mZ_0) = H^0(mZ_0; \operatorname{ad}V\otimes \scrO_X(mZ_0)|mZ_0).$$
Using the exact sequence
$$0 \to H((m-1)Z_0) \to H(mZ_0) \to H^0(Z_0;\operatorname{ad}V_0\otimes
\scrO_X(mZ_0)|Z_0),$$
we see that it suffices to show that, for a generic $V$,
$H^0(Z_0;\operatorname{ad}V\otimes
\scrO_X(mZ_0)|Z_0) =0$ for all $m\geq 1$. For a fixed $m$, the condition that
$H^0(Z_0;\operatorname{ad}V\otimes
\scrO_X(mZ_0)|Z_0) \neq 0$ is a closed condition. Thus since the moduli space
cannot be a countable union of proper subvarieties, it will suffice to show
that the set of $V$ for which $H^0(Z_0;\operatorname{ad}V\otimes
\scrO_X(mZ_0)|Z_0) = 0$ is nonempty for every $m$. Let $\Cal S$ be the germ of
the versal deformation of $V_0|Z_0$ keeping $\det V_0|Z_0$ fixed. By the
assumption that the map from the germ of the versal deformation of $V_0$ to
that
of $V_0|Z_0$ is submersive, it will suffice to show that, for each $m\geq 1$,
the set of
$W\in \Cal S$ such that $H^0(Z_0;\operatorname{ad}W\otimes
\scrO_{Z_0}(mZ_0)) = 0$ is nonempty. One natural method for doing so is to
exhibit a deformation from $V_0|Z_0$ to the $W$ constructed in the course of
Lemma 2.7; roughly speaking this amounts to the claim that the ``moduli space"
of vector bundles on the scheme $Z_0$ is connected. Although we shall proceed
slightly differently, this is the main idea of the argument.

Choose an ample line bundle $\lambda$ on $Z_0$. After passing to some power, we
may assume that both $(V_0|Z_0)\otimes \lambda$ and $W\otimes \lambda$ are
generated by their global sections. A standard argument shows that, in this
case, both $V_0|Z_0$ and $W$ can be written as an extension of $\mu \otimes
\lambda$ by $\lambda  ^{-1}$: Working with $W$ for example, we must show that
there is a map $\lambda ^{-1}\to W$, corresponding to a section of $W\otimes
\lambda$, such that the quotient is again a line bundle. It suffices to show
that there exists a section $s\in H^0(Z_0; W\otimes \lambda)$ such that, for
each $z\in Z_0$, $s(z)\neq 0$ in the fiber of $W\otimes \lambda$ over $z$. Now
for $z$ fixed, the set of $s \in H^0(Z_0; W\otimes \lambda)$ such that $s(z)
=0$
has codimension two in $H^0(Z_0; W\otimes \lambda)$ since $W\otimes \lambda$ is
generated by its global sections. Thus the set of $s \in H^0(Z_0; W\otimes
\lambda)$ such that $s(z) =0$ for some $z\in Z_0$ has codimension at least one,
and so there exists an $s$ as claimed.

Now let $W_0 = \lambda  ^{-1} \oplus (\mu \otimes \lambda)$. Let $(\Cal S_0,
s_0)$ be the germ of the versal deformation of $W_0$ (with fixed determinant
$\mu$). As
$Z_0$ has dimension one, $\Cal S_0$ is smooth. Both $V_0|Z_0$ and $W$
correspond
to extension classes $\xi, \xi ' \in \operatorname{Ext}^1(\mu \otimes \lambda,
\lambda ^{-1})$. Replacing, say, $\xi$ by the class $t\xi, t\in \Cee ^*$, gives
an isomorphic bundle. In this way we obtain a family of bundles $\Cal V$ over
$Z_0 \times \Cee$, such that the restriction of $\Cal V$ to $Z_0 \times t$ is
$V_0|Z_0$ if $t\neq 0$ and is $W_0$ if $t=0$. Hence in the germ $\Cal S_0$
there is a subvariety containing $s_0$ in its closure and consisting of bundles
isomorphic to
$V_0|Z_0$, and similarly for $W$. As $H^0(Z_0;\operatorname{ad}W\otimes
\scrO_{Z_0}(mZ_0)) = 0$, the locus of bundles $U$ in $\Cal S_0$ for which
$H^0(Z_0;\operatorname{ad}U\otimes
\scrO_{Z_0}(mZ_0)) = 0$ is a dense open subset. Since $\Cal S_0$ is a smooth
germ, it follows that there is a small deformation of $V_0|Z_0$ to such a
bundle. Thus the generic small deformation $U$ of $V_0|Z_0$ satisfies $H^0(Z_0;
\operatorname{ad}U\otimes \scrO_{Z_0}(mZ_0)) = 0$, and so the generic $V\in
\frak M_L(w,p)$ has the property that
$H^0(Z_0;\operatorname{ad}V\otimes
\scrO_X(mZ_0)|Z_0) = 0$ for all $m\geq 1$ as well. As we saw above, this
implies the vanishing of $H^0(mZ_0; \operatorname{ad}V\otimes
\scrO_X(mZ_0)|mZ_0)$.
\endproof

\noindent {\bf Remark.} (1) Suppose that $\bar X$ is a singular surface, but
that $\varphi \: X \to \bar X$ is not the minimal resolution. We may still
define the fundamental cycle $Z_0$ for the resolution $\varphi$. Moreover it is
easy to see that $Z_0 \cdot E = 0$ for every component of a generalized
exceptional curve contained in $\varphi ^{-1}(x)$. Thus the hypothesis of (ii)
of Theorem 2.5 implies that $w\cdot E=0$ for such curves.

\smallskip

(2) We have only considered contractions of a very special type, and have
primarily been interested in the case where $w$ is the mod two reduction of
$[K_X]$. However it is natural to ask if the analogues of Theorem 2.3 and 2.5
(and thus Theorem 1.4) holds for more general contractions and choices of $w$,
provided of course that no smooth rational curve of self-intersection $-1$ is
contracted to a point. Clearly the proof of Theorem 2.5 applies to a much wider
class of singularities. Indeed a little work shows that the proof goes over
(with some modifications in case there are components of arithmetic genus one)
to handle the case where we need only assume condition (ii) of (2.5) for those
components $D$ which are smooth rational curves. Another case where it is easy
to check that the conclusions of (2.5) hold is where $w$ is arbitrary and the
dual graph of the singularity is of type $A_k$. We make the following rather
natural conjecture:

\medskip
\noindent {\bf Conjecture 2.9.} The conclusions of Theorem \rom{1.4} hold for
arbitrary choices of $w$ and $\varphi$, provided that $\varphi$ does not
contract any exceptional curves of the first kind.

\section{3. Nonexistence of embedded $\boldkey 2$-spheres}

\ssection{3.1. A base point free theorem}

\theorem{3.1} Let $\pi \: X \to  X'$ be a birational morphism from
the smooth surface $X$ to a normal surface $X'$, not necessarily projective.
Suppose that
$X$ is a minimal surface of  general type, and that $p\in X'$ is
an isolated singular point which is a nonrational singularity. Let $\pi
^{-1}(p) = \bigcup _i D_i$. Then:
\roster
\item"{(i)}"  There exist nonnegative integers $n_i$ with
$n_i >0$ for at least one $i$ such that $K_X+\sum _in_i D_i$ is nef and big.
\item"{(ii)}" Suppose that $q(X) =0$. Then there further exists a choice of $D=
\sum _in_iD_i$ satisfying
\rom{(i)} with $D$ connected and such that there exists a section of
$K_X + D$ which is nowhere vanishing in a neighborhood of
$$E=\bigcup\{\,D_j: (K_X + D)\cdot D_j = 0\,\}.$$
In this case either $E= \emptyset$ or $E=\operatorname{Supp}D$ and $D$ is the
fundamental cycle of the minimal resolution of a minimally elliptic
singularity.
\item"{(iii)}" With $D$ satisfying \rom{(i)} and \rom{(ii)}, the linear system
$K_X+D$ is eventually base point free. Moreover, if $\varphi \: X \to \bar X$
is the associated contraction, then $\bar X$ is a normal projective surface all
of whose singular points are either rational or minimally elliptic.
\endroster
\endproclaim
\proof To prove (i), consider the set of all effective cycles $D = \sum
_ia_iD_i$, where the
$a_i$ are nonnegative integers, not all zero, and such that $h^1(\scrO_D) \neq
0$. This set is not empty by the definition of a nonrational singularity, and
is partially ordered by $\leq$, where $D'\leq D$ if $D-D'$ is effective. Choose
a minimal element $D$ in the set. This means that $D = \sum _in_iD_i$ where
either $n_i = 1$ for exactly one $i$ and $h^1(\scrO_{D_i}) \neq 0$, or for
every irreducible $D_i$ contained in the support of $D$, $D-D_i=D'$ is
effective and $h^1(\scrO_{D'})= 0$. If $D''$ is then any nonzero effective
cycle with
$D''< D$, then there exists an $i$ such that $\scrO_{D-D_i} \to \scrO_{D''}$ is
surjective. By a standard argument, $H^1(\scrO_{D-D_i}) \to H^1(\scrO_{D''})$
is surjective and thus
$h^1(\scrO_{D''})=0$ for every nonzero effective $D''<D$. Finally note that $D$
is connected, since otherwise we could replace $D$ by some connected component
$D_0$ with $h^1(\scrO_{D_0})\neq 0$.

Next we claim that $K_X+D$ is nef. Since $K_X$ is nef, it is clear that
$(K_X+D)\cdot C \geq 0$ for every irreducible curve $C$ not contained in
the support of
$D$, and moreover, for such curves $C$, $(K_X+D)\cdot C = 0$ if and only if $C$
is a smooth rational curve of self-intersection $-2$ disjoint from the support
of $D$. Next suppose that $D_i$ is a curve in the support of $D$ and consider
$(K_X+D)\cdot D_i$.  If
$D= D_i$ then $K_X+D_i |D_i = \omega _{D_i}$, the dualizing sheaf of $D_i$, and
this has degree $2p_a(D_i)-2 \geq 0$ since $p_a(D_i) =h^1(\scrO_{D_i}) > 0$.
Otherwise let
$D' = D-D_i$ consider the exact sequence
$$0 \to \scrO_{D_i}(-D') \to \scrO_D \to \scrO_{D'} \to 0.$$
Thus the natural map $H^1(\scrO_{D_i}(-D')) \to H^1(\scrO_D)$ is surjective
since $H^1(\scrO_{D'})=0$, and so $H^1(\scrO_{D_i}(-D'))\neq 0$ as
$H^1(\scrO_D)\neq 0$. By duality $H^0(D_i; \omega _{D_i}\otimes
\scrO_{D_i}(D'))\neq 0$. On the other hand $\omega _{D_i} = K_X+D_i|D_i$, and
so $\deg (K_X+D'+D_i)|D_i = (K_X+D)\cdot D_i \geq 0$; moreover $(K_X+D)\cdot
D_i
= 0$ only if the divisor class $K_X+D|D_i$ is trivial.

Next, $(K_X+D)^2 \geq K_X^2 >0$, so that $K_X+D$ is big. In fact,
$$(K_X+D)^2 = K_X\cdot (K_X+D)+(K_X+D)\cdot D\geq K_X^2 + K_X\cdot D \geq
K_X^2.$$ Thus $K_X+D$ is big.

To see (ii), let $E = \bigcup\{\,D_j\subseteq \operatorname{Supp}D: (K_X +
D)\cdot D_j = 0\,\}$. We shall also view $E$ as a reduced divisor. We claim
that $\scrO_E(K_X+D) =
\scrO_E$. First assume that $E = D$ (and thus in particular that $D$ is
reduced); in this case we need to show that
$\omega _D = \scrO_D$. By assumption
$D$ is connected. Then $\omega _D$ has degree zero on every reduced irreducible
component of $D$, and by Serre duality $\chi (\omega _D) = - \chi (\scrO_D) =
\frac12(K_X+D)\cdot D=0$. As
$h^1(\omega _D) = h^0(\scrO_D) = 1$,
$h^0(\omega _D) = 1$ as well. As $\omega _D$ has degree zero on every component
of $D$, if $s$ is a section of $\omega _D$, then the restriction of $s$ to
every component $D_i$ of $D$ is either identically zero or nowhere vanishing.
Thus if $s$ is nonzero, since $D$ is connected, $s$ must be nowhere vanishing.
It follows that the map
$\scrO_D\to
\omega _D$ is surjective and is thus an isomorphism.

If $D\neq E$, we apply the argument that showed above that $(K_X+D)\cdot D_i
\geq 0$ to each connected component $E_0$ of the divisor $E$, with $D'=D-E_0$,
to see that there is a section of
$\scrO_{E_0}(K_X+D)$. Since $\scrO_{E_0}(K_X+D)$ has degree zero on each
irreducible component of
$E_0$, the argument that worked for the case $D=E$ also works in this case.

Now let us show that, provided $q(X)=0$, a nowhere zero section of
$\scrO_E(K_X+D) =
\scrO_E$ lifts to a section of $K_X+D$.  It suffices to show that, for every
connected component $E_0$ of $E$, a nowhere vanishing section of
$\scrO_{E_0}(K_X+D)$ lifts to a section of $K_X+D$. Let $D' = D-E_0$.
If $D'=0$ then $D=E=E_0$ and we ask if the map $H^0(\scrO_X(K_X+D)) \to
H^0(\scrO_D(K_X+D))$ is surjective. The cokernel of this map lies in
$H^1(K_X) = 0$ since $X$ is regular. Otherwise $D' \neq 0$. Beginning with the
exact sequence
$$0 \to \scrO_{D'}(-E_0)\to \scrO_D \to \scrO_{E_0} \to 0,$$
and tensoring with $\scrO_X(K_X+D)$, we obtain the exact sequence
$$0 \to \scrO_{D'}(K_X+D-E_0)\to \scrO_D(K_X+D) \to \scrO_{E_0} \to 0.$$
Now $\scrO_{D'}(K_X+D-E_0) = \scrO_{D'}(K_X+D')= \omega _{D'}$ and by duality
$h^0(\omega _{D'}) = h^1(\scrO_{D'}) = 0$. Thus $H^0(\scrO_D(K_X+D))$ includes
into $H^0(\scrO_{E_0}) = \Cee$ and so it suffices to prove that
$H^0(\scrO_D(K_X+D)) \neq 0$, in which case it has dimension one. On the other
hand, using the exact sequence
$$0 \to \scrO_X(K_X) \to \scrO_X(K_X+D) \to \scrO_D(K_X+D) \to 0,$$
we see that $h^0(\scrO_D(K_X+D)) \geq h^0(\scrO_X(K_X+D))-p_g(X)$. Since
$h^2(K_X+D) = h^0(-D)=0$, the Riemann-Roch theorem implies that
$$h^0(\scrO_X(K_X+D)) = h^1(\scrO_X(K_X+D)) + \frac12(K_X+D)\cdot D + 1 +
p_g(X).$$
Since all the terms are positive, we see that indeed
$h^0(\scrO_X(K_X+D))-p_g(X)\geq 1$, and that $h^0(\scrO_X(K_X+D))-p_g(X)=1$ if
and only if $h^1(\scrO_X(K_X+D)) = 0$ and $(K_X+D)\cdot D_i = 0$ for every
component $D_i$ contained in the support of $D$. This last condition says
exactly that $E= \operatorname{Supp}D$, and thus, as $D$ is connected,
that $E_0=E$. We claim that in this last case
$D$ is minimally elliptic. Indeed, for every effective divisor $D'$ with
$0<D'<D$, we have
$$p_a(D') = 1 - h^0(\scrO_{D'}) + h^1(\scrO_{D'})= 1 - h^0(\scrO_{D'})\leq 0.$$
Thus $D$ is the fundamental cycle for the resolution of a minimally elliptic
singularity.

Finally we prove (iii). The irreducible curves $C$ such that $(K_X+D)\cdot C=0$
are the components $D_i$ of the support of $D$ such that $(K_X+D)\cdot D_i =
0$,
as well as smooth rational curves of self-intersection
$-2$ disjoint from $\operatorname{Supp}D$. These last contribute
rational double points, so that we need only study the $D_i$ such that
$(K_X+D)\cdot D_i = 0$. We have seen in (ii) that either there are no such
$D_i$, or every $D_i$ in the support of $D$ satisfies $(K_X+D)\cdot D_i = 0$
and the contraction of $D$ is a minimally elliptic singularity.

Let $\bar X$ be the normal surface obtained by contracting all the
irreducible curves $C$ on
$X$ such that  $(K_X+D)\cdot C=0$. The line bundle $\scrO_X(K_X+D)$ is trivial
in a neighborhood of these curves, either because they correspond to a rational
singularity or because we are in the minimally elliptic case and by (ii). So
$\scrO_X(K_X+D)$ induces a line bundle on $\bar X$ which is ample, by the
Nakai-Moishezon criterion. Thus $|k(K_X+D)|$ is base point free for all
$k\gg0$.
\endproof

\ssection{3.2. Completion of the proof}

We now prove  Theorem 1.5:

\theorem{1.5} Let $X$ be a minimal simply connected algebraic surface of
general type, and let
$E\in H^2(X; \Zee)$ be a $(1,1)$-class satisfying $E^2=-1$, $E\cdot K_X = 1$.
Let $w$ be the \rom{mod} $2$ reduction of $[K_X]$. Then there exist:
\roster
\item"{(i)}"  an integer $p$
and \rom(in case $p_g(X)=0$\rom) a chamber $\Cal C$ of type $(w,p))$ and
\item"{(ii)}" a $(1,1)$-class $M\in H^2(X; \Zee)$
\endroster
such that
$M\cdot E=0$ and $\gamma _{w,p}(X)(M^d) \neq 0$ \rom(or, in case $p_g(X)=0$,
$\gamma _{w,p}(X; \Cal C)(M^d) \neq 0$\rom).
\endstatement
\noindent {\it Proof.}  We begin with the following lemma:

\lemma{3.2} With $X$ and $E$ as above, there exists an orientation preserving
diffeomorphism
$\psi \: X \to X$ such that
$\psi ^*[K_X] = [K_X]$ and such that $\psi ^*E\cdot [C] \geq 0$ for every
smooth rational curve $C$ on $X$ with $C^2 = -2$.
\endstatement
\proof Let $\Delta = \{[C_1], \dots, [C_k]\}$ be the set of smooth
rational curves on $X$ of self-intersection $-2$, and let $r_i\: H^2(X; \Zee)
\to H^2(X; \Zee)$ be the reflection about the class $[C_i]$. Then $r_i$ is
realized by an orientation-preserving self-diffeomorphism of $X$, $r_i^*[K_X] =
[K_X]$, and $r_i$ preserves the image of $\operatorname{Pic}X$ inside $H^2(X;
\Zee)$. Let $\Gamma$ be the finite group generated group by the $r_i$. Since
the classes $[C_i]$ are linearly independent, the set
$$\{\, x\in H^2(X; \Ar)\: x\cdot [C_i] \geq 0\,\}$$
has a nonempty interior. Moreover, if $\Delta ' = \Gamma \cdot \Delta$, and we
set $W^\delta = \delta ^\perp$ for $\delta \in \Delta'$, then the connected
components of the set $H^2(X; \Ar) - \bigcup _{\delta \in \Delta '}W^\delta$
are the fundamental domains for the action of $\Gamma$ on $H^2(X; \Ar)$.
Clearly at least one of these connected components lies inside $\{\, x\in
H^2(X; \Ar)\: x\cdot [C_i] \geq 0\,\}$. Thus given $E$ (or indeed an arbitrary
element of $H^2(X; \Ar)$), there exists a $\gamma \in \Gamma$ such that $\gamma
(E)\cdot [C_i] \geq 0$ for all $i$. As every $\gamma \in \Gamma$ is realized by
an orientation preserving self-diffeomorphism $\psi$, this concludes the
proof of (3.2).
\endproof

Thus, to prove Theorem 1.5, it is sufficient by the naturality of the
Donaldson polynomials to prove it for every class
$E$ satisfying  $E^2=-1$, $E\cdot K_X = 1$, and $E\cdot [C] \geq 0$ for every
smooth rational curve $C$ on $X$ with $C^2 = -2$. We therefore make this
assumption in what follows.
Given Theorem  1.4, it therefore suffices to find a nef and big divisor $M$
orthogonal to $E$, which is eventually base point free, such that the
contraction morphism defined by $|kM|$ has an image with at worst rational and
minimally elliptic singularities (note that, since $X$ is assumed minimal, no
exceptional curves can be contracted). Thus we will be done by the following
lemma:

\lemma{3.3} There exists a nef and big divisor
$M$ which is eventually base point free and such that
\roster
\item $M\cdot E = 0$.
\item The contraction $\bar X$ of $X$ defined by $|kM|$ for all $k\gg 0$ has
only rational and minimally elliptic singularities.
\endroster
\endstatement
\proof To find $M$ we proceed as follows: consider the
divisor
$K_X + E=M$. As
$K_X\cdot E=1$ and $E^2 = -1$, $M$ is orthogonal to $E$. Moreover
$M^2=(K_X+E)^2 = K_X^2 + 1 > 0$. We now consider separately the cases where $M$
is nef and where $M$ is not nef.

\medskip
\noindent {\bf Case I:} $M= K_X+E$ is nef.

Consider the union of all the curves $D$ such that $M \cdot D = 0$. The
intersection matrix of the $D$ is  negative definite, and so we can contract
all the $D$ on $X$ to obtain a normal surface $X'$. If $X'$ has only rational
singularities, then $M$ induces an ample divisor on $X'$ and so $M$ itself is
eventually base point free. In this case we are done. Otherwise we may apply
Theorem 3.1 to find a subset $D_1, \dots, D_t$ of the curves $D$ with $M\cdot
D=0$ and positive integers
$a_i$ such that the divisor
$K_X+\sum _ia_iD_i$ is nef, big, and eventually base point free, and such that
the contraction $\bar X$ of
$X$ has only rational and minimally elliptic singularities, with exactly one
nonrational singularity. Note that
$D_i\cdot E = -D_i \cdot K_X \leq 0$, and
$D_i \cdot E = 0$ if and only if $D_i \cdot K_X = 0$, or in other words if and
only if $D_i$ is a smooth rational curve of self-intersection $-2$. Setting $e
= -\sum _ia_i(D_i \cdot E)$, we have $e \geq 0$, and $e=0$ if and only if
$D_i
\cdot E =0$ for all $i$.  But as $\bar X$ has a nonrational singularity, we
cannot have $D_i
\cdot E =0$ for all $i$, for then all singularities would be rational double
points. Thus $e>0$. Now the $\Bbb Q$-divisor $M'= K_X + \frac1{e}\sum _ia_iD_i$
is a rational convex combination of $K_X$ and $K_X+\sum _ia_iD_i$, and $M'
\cdot E=0$. Moreover either $M'$ is a strict convex combination of $K_X$ and
$K_X+\sum _ia_iD_i$ (if $e>1$) or $M'=K_X+\sum _ia_iD_i$ (if $e=1$). In the
second case,
$M'$ satisfies (1) of Lemma 3.3, and it is eventually base point free by (iii)
of Theorem 3.1. Thus $M'$ satisfies the conclusions of Lemma 3.3. In the first
case,
$M'$ is nef and big, and the only curves $C$ such that $M' \cdot C=0$ are
curves $C$ such that $K_X
\cdot C = 0$ and $(K_X+\sum _ia_iD_i )\cdot C=0$. The set of all such
curves must therefore be a subset of the set of all  smooth rational curves on
$X$ with self-intersection $-2$. Hence, if $X''$ denotes the contraction of all
the curves $C$ on $X$ such that $M'\cdot C=0$, then $X''$ has only rational
singularities and $M'$ induces an ample $\Bbb Q$-divisor on $X''$. Once again
some multiple of
$M'$ is eventually base point free and (1) and (2) of Lemma 3.3 are satisfied.
Thus we have proved the lemma in case $K_X+E$ is nef.

\medskip
\noindent {\bf Case II:} $M= K_X+E$ is not nef.

Let $D$ be an irreducible curve with $M\cdot D<0$. We claim first that in this
case $D^2<0$. Indeed, suppose that $D^2 \geq 0$. As $\operatorname{Pic}X\otimes
_\Zee \Ar$ has signature $(1, \rho -1)$, the set
$$\Cal Q = \{\, x\in \operatorname{Pic}X\otimes _\Zee \Ar\: x^2 \geq 0, x\neq
0\,\}$$ has two connected components, and two classes $x$ and $x'$ are in the
same connected component of $\Cal Q$ if and only if $x\cdot x' \geq 0$ (cf\.
[13] p\. 320 Lemma 1.1). Now $(K_X+E)\cdot K_X = (K_X+E)^2 = K_X^2 +1 >0$, so
that
$K_X+E$ and
$K_X$ lie in the same connected component of $\Cal Q$. Likewise, if $D^2 \geq
0$, then since $K_X\cdot D \geq 0$, $K_X$ and $D$ lie in the same connected
component of $\Cal Q$. Thus $D$ and $K_X+E$ lie in the same connected component
of $\Cal Q$, so that $(K_X+E)\cdot D \geq 0$. Conversely, if $M\cdot D<0$, then
$D^2<0$.

Fix an irreducible curve $D$ with $M\cdot D<0$, and let $d = -E\cdot D
>K_X\cdot D \geq 0$. Recall that by assumption $E\cdot D \geq 0$ if $D$ is a
smooth rational curve of self-intersection $-2$. If
$p_a(D) \geq 1$, then set $M' = K_X+ \frac1{d}D$. Then $M'\cdot E = 0$ by
construction. Moreover we claim that $M'$ is nef and big. Indeed
$$(M')^2 = (K_X+ \frac1{d}D)^2 = K_X\cdot (K_X+ \frac1{d}D) + \frac1{d}(K_X+
\frac1{d}D)\cdot D.$$
Thus $M'$ is big if it is nef and to see that
$M'$ is nef it suffices to show that $M'\cdot D \geq 0$. But
$$M'\cdot D = K_X\cdot D + \frac1{d}D^2 = 2p_a(D) - 2 - \left(1-
\frac1{d}\right)D^2.$$
As $D^2 < 0$, we see that $M'\cdot D\geq 0$, and $M'\cdot D=0$ if and only if
$p_a(D) = 1$ and $d=1$. Suppose that $p_a(D) = 1$ and $M'\cdot D=0$.
Using the exact sequence
$$0 \to \scrO_X(K_X) \to \scrO_X(K_X+D)\to \omega _D \to 0,$$
and arguments as in the proof of Theorem 3.1, we see that the linear system
$M'$ is eventually base point free and that the associated contraction has just
rational double points and a minimally elliptic singular point which is the
image of $D$. In all other cases,
$M'\cdot D >0$, so that the curves orthogonal to $M'$ are smooth rational
curves of self-intersection
$-2$. Again, some positive multiple of $M'$ is eventually base point free and
the contraction has just rational singularities.

Thus we may assume that $p_a(D) = 0$ for every irreducible curve $D$ such that
$M\cdot D <0$. By assumption $D^2 \neq -1, -2$, so that $D^2 \leq -3$. Thus $d=
- -D\cdot E \geq 2$. If either $D^2 \leq -4$ or $D^2 = -3$ and $d \geq 3$, then
again let $M' = K_X+ \frac1{d}D$. Thus $M'\cdot E =0$ and
$$M'\cdot D = K_X\cdot D + \frac1{d}D^2 = - 2 - \left(1-
\frac1{d}\right)D^2\geq 0.$$
Thus $M'$ is nef and big, and some multiple of $M'$ is eventually base point
free, and the associated contraction has just rational singularities. The
remaining case is where there is a smooth rational curve $D$ on $X$ with
self-intersection
$-3$ and such that
$-D\cdot E = 2$. In this case $K_X \cdot D = 1$, and so $D-E$ is orthogonal to
$K_X$. Note that $D-E$ is not numerically trivial since $D$ is not numerically
equivalent to $E$. Thus, by the Hodge index theorem
$(D-E)^2 <0$. But
$$(D-E)^2 = -3 + 4 -1 =0,$$
a contradiction. Thus this last case does not arise.
\endproof

\section{Appendix: On the canonical class of a rational surface}

Let $\Lambda _n$ be a lattice of type $(1,n)$, i.e. a free $\Bbb Z$-module of
rank $n+1$, together with a quadratic form $q\: \Lambda _n\to \Bbb Z$, such
that there exists an orthogonal basis $\{ e_0, e_1, \dots , e_n\}$ of $\Lambda
_n$ with
$q(e_0) = 1$ and $q(e_i) = -1$ for all $i>0$. Fix once and for all such a
basis. We shall always view $\Lambda _n$ as included in $\Lambda _{n+1}$ in
the obvious way. Let  $$\kappa _n = 3e_0 - \sum _{i=1}^ne_i.$$ Then $q(\kappa
_n) = 9-n$ and
$\kappa _n$ is characteristic, i.e. $\kappa _n\cdot\alpha  \equiv q(\alpha)
\mod 2$ for all $\alpha \in \Lambda _n$.

The goal of this appendix is to give a proof, due to the first author, R\.
Miranda, and J\.W\. Morgan, of the following:

\theorem{A.1} Suppose that $n\leq 8$ and that $\kappa\in \Lambda _n$ is a
characteristic vector satisfying $q(\kappa ) = 9-n$. Then there exists an
automorphism $\varphi$ of $\Lambda _n$ such that $\varphi (\kappa) = \kappa
_n$.
A similar statement holds for $n=9$ provided that $\kappa$ is primitive.
\endproclaim
\demo{Proof} We shall freely use the notation and results of Chapter II of
[13] and shall quote the results there by number. For the purposes of the
appendix, chamber shall mean a chamber in $\{\,x\in \Lambda _n \otimes \Bbb R
\mid x^2 =1\,\}$ for the set of walls defined by the set $\{\, \alpha \in
\Lambda _n \mid \alpha ^2 = -1\,\}$. Let
$C_n$ be the chamber associated to $\kappa _n$ [13, p\. 329, 2.7(a)]: the
oriented walls of $C_n$ are exactly the set
$$\{\,\alpha\in \Lambda _n\mid q(\alpha ) = -1,\alpha\cdot \kappa _n =1\,\}.$$
Then
$\kappa _n$ lies in the interior of
$\Bbb R^+\cdot C_n$, by [13, p\. 329, 2.7(a)]. Similarly $\kappa$ lies in the
interior of a set of the form  $\Bbb R^+\cdot C$ for some chamber $C$, since
$\kappa$ is not orthogonal to any wall (because it is characteristic) and
$q(\kappa )> 0$. But the automorphism group of $\Lambda _n$ acts transitively
on the chambers, by [13 p. 324]. Hence we may assume that $\kappa \in C_n$. In
this case we shall prove that $\kappa =
\kappa _n$. We shall refer to $C_n$ as the {\sl fundamental chamber} of
$\Lambda _n\otimes _{\Bbb Z}\Bbb R$. Let us record two lemmas about $C_n$.

\lemma{A.2} An automorphism $\varphi$ of $\Lambda _n$ fixes $C_n$ if and
only if it fixes $\kappa _n$.
\endproclaim
\demo{Proof} The oriented walls of $C_n$ are precisely
the
$\alpha \in \Lambda _n$ such that $q(\alpha) = -1$ and $\kappa _n\cdot \alpha
=1$. Thus, an automorphism fixing $\kappa _n$ fixes $C_n$. The converse follows
from [13, p\. 335, 4.4].
\endproof

\lemma{A.3} Let $\alpha = \sum _i\alpha _ie_i$ be an oriented wall of
$C_n$, where $e_0, \dots , e_n$ is the standard basis of $\Lambda _n$. After
reordering the elements $e_1, \dots, e_n$, let us assume that
$$|\alpha _1| \geq |\alpha _2| \geq \dots \geq |\alpha _n|.$$
Then for $n\leq 8$,
the possibilities for $(\alpha _0, \dots \alpha _n)$ are as follows
\rom(where we omit the $\alpha _i$ which are zero\rom):
\roster
\item $\alpha _0 = 0, \alpha _1 = 1$;
\item $\alpha _0 =1, \alpha _1 = \alpha _2  = -1$ $(n\geq 2)$;
\item $\alpha _0 = 2, \alpha _1 = \alpha _2  = \alpha _3 = \alpha _4  = \alpha
_5 =-1$ $(n \geq 5)$;
\item $\alpha _0 = 3, \alpha _1 = -2, \alpha _2  = \alpha _3 = \alpha _4  =
\alpha _5 =\alpha _6 = \alpha _7=-1$ $(n \geq 7)$;
\item $\alpha _0 = 4, \alpha _1 = \alpha _2  = \alpha _3 =  2, \alpha _4 =
\dots = \alpha _8 = -1$  $(n=8)$;
\item $\alpha _0 = 5, \alpha _1 = \dots = \alpha _6 = -2, \alpha _7 =
\alpha _8 = -1$ $(n=8)$;
\item $\alpha _0 =  6, \alpha _1 = -3, \alpha _2 = \dots = \alpha _8 = -2$
$(n=8)$.
\endroster
\endproclaim
\demo{Proof} This statement is extremely well-known as the characterization of
the lines on a del Pezzo surface (see [7], Table 3). We can
give a proof as follows. It clearly suffices to prove the result for $n=8$. But
for $n=8$, there is a bijection between the $\alpha$ defining an oriented wall
of $C_8$ and the elements $\gamma \in \kappa _8^{\perp}$ with $q(\gamma )=-2$.
This bijection is given as follows: $\alpha$ defines an oriented wall
of $C_8$ if and only if $q(\alpha ) = -1$ and $\kappa _8\cdot \alpha = 1$.
Map $\alpha$ to $\alpha - \kappa _8 = \gamma$. Thus, as $q(\kappa _8 ) = 1$,
$q(\gamma ) = -2$ and $\gamma \cdot \kappa _8 = 0$. Conversely, if $\gamma \in
\kappa _8^{\perp}$ satisfies $q(\gamma )=-2$, then $\gamma + \kappa _8$ defines
an oriented wall of $C_8$.

Now the number of $\alpha$ listed above, after we are allowed to reorder the
$e_i$, is easily seen to be
$$8+\binom 82 +\binom 85 +8\cdot 7 +\binom 83 +\binom 82 + 8 = 240.$$
Since this is exactly the number of vectors of square $-2$ in $-E_8$, by e\.g\.
[36], we must have enumerated all the possible $\alpha$.
\endproof

Write $\kappa = \sum _{i=0}^na_ie_i$, where $e_i$ is the standard basis of
$\Lambda _n$ given above. Since $\kappa \cdot e_i>0$, $a_i <0$. After
reordering the elements $e_1, \dots, e_n$, we may assume that
$$|a_1|\geq |a_2| \geq \dots \geq |a_n|.$$
By inspecting  the cases in Lemma
A.3, for every
$\alpha = \sum _i\alpha
_ie_i$ not of the form $e_i$,
$\alpha _i \leq 0$ for all $i\geq 1$. Given $\alpha = \sum _i\alpha
_ie_i$ with $\alpha
\neq e_i$ for any
$i$, let  us call  $\alpha$ {\sl
well-ordered} if
$$|\alpha_1|\geq |\alpha_2| \geq \dots \geq |\alpha_n|.$$
Quite generally, given $\alpha = \alpha _0e_0 + \sum _{i>0}\alpha
_ie_i$, we define the {\sl reordering} $r(\alpha)$ of $\alpha$ to be
$$r(\alpha ) = \alpha _0e_0 + \sum
_{i>0}\alpha _{\sigma (i)}e_i,$$
where $\sigma $ is a permutation of $\{ 1,
\dots , n\}$ such that $r(\alpha)$ is well-ordered. Clearly $r(\alpha)$ is
independent of the choice of $\sigma$.

We then have the following:

\claim{A.4} $\kappa \in C_n$ if and only if $\kappa \cdot \alpha >0$ for every
well-ordered wall $\alpha$.
\endstatement
\proof Clearly if $\kappa \in C_n$, then $\kappa \cdot \alpha >0$ for
every $\alpha$, well-ordered or not. Conversely, suppose that $\kappa \cdot
\alpha >0$ for every well-ordered wall $\alpha$.
We claim that
$$\alpha \cdot \kappa \geq r(\alpha )\cdot \kappa,\tag{$*$}$$
which clearly implies (A.4) since $r(\alpha)$ is well-ordered.
Now
$$\alpha \cdot \kappa = \alpha _0 a_0 - \sum _{i>0}\alpha _i a_i.$$
Since $\alpha _i < 0$ and $a_i <0$,
($*$) is easily reduced to the following statement about positive real numbers:
if $c_1\geq \dots \geq c_n$ is a sequence of positive real numbers and $d_1,
\dots, d_n$ is any sequence of positive real numbers, then a permutation
$\sigma$ of  $\{1,
\dots , n\}$ is such that $\sum _ic_id_{\sigma (i)}$ is maximal exactly when
$d_{\sigma (1)} \geq \dots \geq d_{\sigma (n)}$. We leave the proof of this
elementary fact to the reader.
\endproof

Next, we claim the following:

\lemma{A.5} View $\Lambda _n\subset \Lambda _{n+1}$.
Defining $\kappa _{n+1}$ and $C_{n+1}$ in the natural way for $\Lambda _{n+1}$,
suppose that $\kappa \in C_n$. Then $\kappa ' =\kappa - e_{n+1} \in C_{n+1}$.
\endproclaim
\demo{Proof} We have ordered our basis $\{e_0, \dots, e_n\}$ so that
$$|a_1|\geq |a_2| \geq \dots \geq |a_n|.$$
Since $a_i <0$ for all $i$, $|a_n|\geq 1$. Thus the coefficients of $\kappa '$
are also so ordered. Note also that all coefficients of $\kappa
'$ are less than zero, so that the inequalities from (1) of Lemma A.3 are
automatic. Given any other wall $\alpha '$ of $C_{n+1}$, to verify that $\kappa
'\cdot \alpha '> 0$, it suffices to look at $\kappa '\cdot r(\alpha ')$, where
$r(\alpha ')$ is the reordering of $\alpha '$. Expressing $\alpha '$ as a
linear combination of the standard basis vectors, if some coefficient is zero,
then $r(\alpha ') \in \Lambda _n$. Clearly, in this case, viewing $r(\alpha ')$
as an element of $\Lambda _n$, it is a wall of $C_n$. Since then $\kappa '
\cdot r(\alpha ') = \kappa \cdot r(\alpha ')$, we have $\kappa ' \cdot r(\alpha
') >0$ in this case.

In the remaining case, $r(\alpha ')$ does not lie in $\Lambda _n$. This can
only happen for $n=1,4,6,7$, with $\alpha '$ one of the new types of walls
corresponding
to the cases (2) --- (7) of Lemma A.3. Thus, the only thing we need to check is
that, every time we introduce a new type of wall, we still get the inequalities
as needed. Since $r(\alpha ')$ is well-ordered, we can assume that it is in
fact one of the walls listed in  Lemma A.3.

The $n=1$ case simply says that $a_0 > -a_1+1$. However, we can easily solve
the equations $a_0^2 - a_1^2 = 8$, $a_1<0$ to get $a_0 = 3, a_1 =
- -1$. Since $3> 1+1$, we are done in this case.

Next assume
that $n =4$.  We have $\kappa = a_0e_0 +
\sum _{i=1}^4a_ie_i$. We must show that $2a_0> -\sum _{i=1}^4a_i + 1$.
We know that
$a_0>-a_1-a_2$, hence that
$a_0 \geq -a_1-a_2 +1$. Moreover
$a_0\geq -a_3-a_4 +1$ since $|a_1|\geq |a_2| \geq |a_3|  \geq |a_4|$.
Adding gives
$2a_0 \geq -\sum _{i=1}^4a_i + 2 = (-\sum _{i=1}^4a_i+ 1) +1$
and therefore
$2a_0> -\sum _{i=1}^4a_i + 1$.
The case where $n=6$ is similar: we must show that
$3a_0 > -2a_1 -\sum _{i=2}^6 a_i +1$.
But we know that
$2a_0 \geq -\sum _{i=1}^5a_i +1$
and that
$a_0 \geq -a_1 -a_6 +1$.
Adding gives the desired inequality.
For $n=7$, we have three new inequalities to check. The inequality
$$4a_0> -2a_1 - 2a_2 - 2a_3 - \sum _{i=4}^7a_i +1$$
follows by adding the inequalities
$3a_0 > -2a_1 -\sum _{i=2}^7 a_i$
and $a_0> -a_2 - a_3$.
The inequality
$$5a_0 > -2\sum _{i=1}^6a_i -a_7 +1$$
follows from adding
the inequalities
$3a_0 > -2a_1 - \sum _{i=2}^7a_i$
and $2a_0 > -\sum _{i=2}^5a_i$.
Likewise, the last inequality
$$6a_0 > -3a_1 - 2\sum _{i=2}^7a_i+2$$
follows by adding up the three inequalities
$3a_0 > -2a_1 -\sum _{i=2}^7 a_i$,
$2a_0 > -\sum _{i=1}^5a_i$,
and
$a_0 > -a_6-a_7$.
Thus we have established the
lemma.
\endproof

\noindent {\it Completion of the proof of Theorem \rom{A.1}.} Begin with
$\kappa$. Applying Lemma A.5 and induction, if $n<8$, then the vector
$\eta = \kappa -
\sum _{j=n+1}^8e_j$ lies in the fundamental chamber of $\Lambda _8$. Moreover
$\eta$ is a characteristic vector of square $1$. Thus $\eta ^{\perp} \cong
- -E_8$. The same is true for $\kappa _8 = 3e_0 - \sum _{i=1}^8e_i= \kappa _n -
\sum _{j=n+1}^8e_j$. Clearly, then, there is an automorphism $\varphi$ of
$\Lambda _8$ such that
$\eta=\varphi (\kappa _8)$. But both $\eta $ and $\kappa _8$ lie in the
fundamental chamber for
$\Lambda _8$. Since the automorphism group preserves the chamber structure, the
automorphism
$\varphi$ must stabilize the fundamental chamber. By Lemma A.2, $\varphi
(\kappa _8) = \kappa _8$. Thus $\eta =
\kappa _8$. Hence
$\kappa - \sum _{j=n+1}^8e_j = \kappa _n- \sum _{j=n+1}^8e_j$.
It follows that $\kappa = \kappa _n$.
\endproof

\remark{Note} To handle the case $n=9$, we argue that every vector $\kappa\in
\Lambda _9$ which is primitive of square zero and characteristic is conjugate
to $\kappa _9$ as above. To do this, an easy argument shows that, if $\kappa$
is such a class, then there is an orthogonal splitting
$$\Lambda _9 \cong \langle \kappa , \delta\rangle \oplus (-E_8),$$
where $\delta$ is an element of $\Lambda _9$ satisfying $\delta \cdot \kappa =
1$ and $q(\delta )=1$. Thus clearly every two such $\kappa$ are conjugate.
\endremark


\Refs



\ref \no 1\by M. Artin \paper On isolated rational singularities of
surfaces\jour  Amer. J. Math.
\vol 88 \pages 129--136 \yr 1966
\endref

\ref \no 2\by R. Barlow\paper A simply connected surface of general type with
$p_g = 0$\jour Invent. Math. \vol 79\pages 293--301 \yr 1985\endref

\ref \no 3\by W. Barth, C. Peters, A. Van de Ven \book Compact Complex
Surfaces , {\rm Ergebnisse der Mathematik und ihrer Grenz\-gebiete 3.
Folge} {\bf 4}\publ Springer\publaddr Berlin Heidelberg New York Tokyo\yr 1984
\endref

\ref \no  4\by F. Bogomolov \paper Holomorphic tensors and vector bundles
on projective varieties \jour Math. USSR Izvestiya \vol 13 \yr 1979 \pages
499--555\endref

\ref \no 5\by R. Brussee \paper Stable bundles on blown up surfaces \jour Math.
Zeit. \vol 205 \yr 1990 \pages 551--565
\endref

\ref \no 6\bysame \paper On the $(-1)$-curve conjecture of Friedman and
Morgan \jour Invent. Math. \vol 114 \pages 219--229 \yr 1993 \endref

\ref \no 7\by M. Demazure \paper Surfaces de del Pezzo II \inbook
S\'eminaire sur les Singularit\'es des Surfaces {\rm Palaiseau, France
1976--1977} \eds M. Demazure, H. Pinkham, B. Teissier \bookinfo Lecture Notes
in Mathematics {\bf 777}
\publ Springer \publaddr Berlin Heidelberg New York \yr 1980\endref

\ref \no 8\by S.K. Donaldson \paper Anti-self-dual Yang-Mills connections over
complex algebraic surfaces and stable vector bundles\jour Proc. Lond. Math.
Soc. \vol 50\pages 1--26 \yr 1985\endref

\ref \no 9\bysame \paper Irrationality and the $h$-cobordism
conjecture\jour J. Differ. Geom. \vol 26\pages 141--168 \yr 1987\endref

\ref \no 10\bysame\paper Polynomial invariants for smooth
four-manifolds \jour Topology \vol 29 \pages  257--315 \yr 1990\endref

\ref \no 11\by S.K. Donaldson, P. Kronheimer \book The Geometry of
Four-Manifolds
\publ Clarendon Press \publaddr Oxford \yr 1990 \endref

\ref \no 12\by R. Friedman \book Stable Vector Bundles over Algebraic Varieties
\toappear
\endref

\ref \no 13\by R. Friedman, J. W.  Morgan\paper On the
diffeomorphism types of certain algebraic surfaces I \jour J. Differ. Geom.
\vol   27
\pages  297--369 \yr 1988\endref

\ref \no 14\bysame\paper Algebraic surfaces and
$4$-manifolds:  some conjectures and speculations \jour Bull. Amer. Math. Soc.
(N.S.)
\vol   18 \pages  1--19 \yr 1988\endref

\ref \no 15\bysame\book Smooth Four-Manifolds and
Complex Surfaces, {\rm Ergebnisse der Mathematik und ihrer Grenz\-gebiete 3.
Folge} {\bf 27} \publ Springer \publaddr Berlin Heidelberg  New York \yr
1994\endref

\ref \no 16\by Y. Kawamata, K.  Matsuda, K. Matsuki \paper Introduction to the
minimal model problem  \inbook Algebraic Geometry,
{\rm Sendai   Adv. Stud. Pure Math. }\ed T. Oda\vol   10
\pages  283--360\publ Kinokuniya and  Amsterdam North-Holland \publaddr Tokyo
\yr 1987\endref

\ref \no 17\by M. Kneser\paper Klassenzahlen indefiniter quadratischer
Formen in drei oder mehr Ver\"anderlichen \jour Arch. der Math. \vol   7
\pages  323--332
\yr 1956\endref

\ref \no 18\by D. Kotschick \paper On manifolds homeomorphic to $\Bbb CP ^2
\# 8\overline {\Bbb CP}  ^2$\jour  Invent. Math. \vol   95 \pages  591--600
\yr 1989\endref

\ref \no 19\bysame\paper $SO(3)$-invariants for $4$-manifolds
with
$b_2^+ = 1$\jour Proc. Lond. Math. Soc. \vol   63 \pages  426--448 \yr
1991\endref

\ref \no 20\bysame \paper Positivity versus rationality of algebraic
surfaces \toappear \endref

\ref \no 21\by D. Kotschick, J.W. Morgan\paper $SO(3)$-invariants for
$4$-manifolds with $b_2^+ = 1$ II \jour J. Differ. Geom.\vol 39 \yr 1994 \pages
443--456
\endref

\ref \no 22\by H.B. Laufer\paper On minimally elliptic singularities \jour
Amer. J. Math.
\vol   99 \pages  1257--1295 \yr 1977\endref

\ref \no 23\by J. Li \paper Algebraic geometric interpretation of
Donaldson's polynomial invariants \jour J. Differ. Geom. \vol   37 \pages
417--466
\yr 1993\endref

\ref \no 24\bysame \paper Kodaira dimension of moduli space of vector
bundles on surfaces\jour Invent. Math. \vol   115 \pages  1--40 \yr 1994\endref

\ref \no 25\by J.W. Morgan \paper Comparison of the Donaldson polynomial
invariants  with their alge\-bro-geo\-metric analogues \jour Topology \vol   32
\pages  449--488 \yr 1993\endref

\ref \no 26\by J.W. Morgan, T. Mrowka, D. Ruberman \book
The $L^2$-Moduli Space and  a Vanishing Theorem for Donaldson Polynomial
Invariants \toappear \endref

\ref \no 27\by K. O'Grady \paper Moduli of bundles on surfaces: some global
results \toappear\endref

\ref \no 28\by C. Okonek, A. Van de Ven\paper Stable bundles
and differentiable  structures on certain elliptic surfaces\jour  Invent. Math.
\vol   86  \yr 1986\pages  357--370 \endref

\ref \no 29\bysame \paper $\Gamma$-type-invariants associated to
$PU ( 2)$-bundles and the differentiable structure of Barlow's
surface\jour   Invent. Math. \vol   95 \pages  601--614 \yr 1989\endref

\ref \no 30\by V.Y. Pidstrigach\paper Deformation of instanton surfaces\jour
Math. USSR Izvestiya \vol   38 \pages  313--331 \yr 1992\endref

\ref \no 31\by V.Y. Pidstrigach, A.N.  Tyurin \paper Invariants
of the smooth  structure of an algebraic surface arising from the Dirac
operator
\jour  Russian Academy of Science Izvestiya Mathematics, Translations of the
AMS
\vol   40 \pages  267--351 \yr 1993\endref

\ref \no 32\by Z.B. Qin\paper Birational properties of moduli spaces of
stable locally free rank-$2$ sheaves on algebraic surfaces\jour
Manuscripta Math. \vol   72 \pages  163--180 \yr 1991\endref

\ref \no 33\bysame \paper Complex structures on certain differentiable
$4$-manifolds\jour Topology \vol   32 \pages  551--566 \yr 1993\endref

\ref \no 34\bysame \paper On smooth structures of potential surfaces of
general type homeomorphic to rational surfaces\jour
Invent. Math. \vol   113 \pages  163--175 \yr 1993\endref

\ref \no 35\by I. Reider \paper Vector bundles of rank $2$ and linear systems
on algebraic surfaces \jour Annals of Math. \vol 127 \pages 309--316 \yr
1988 \endref

\ref \no 36\by J.P. Serre \book Cours d'Arithm\'etique \publ Presses
Universitaires de France \publaddr Paris \yr 1970 \endref

\ref \no 37\by A. Van de Ven \paper On the differentiable structure of
certain algebraic surfaces\paperinfo  S\'em. Bourbaki vol. 1985--1986
Expos\'es 651--668 n${}^\circ$ {\bf 667} Juin 1986 \jour Ast\'erisque \vol
145--146 \yr 1987\endref

\ref \no 38\by C.T.C. Wall\paper Diffeomorphisms of $4$-manifolds\jour J.
Lond. Math.  Soc. \vol   39 \pages  131--140 \yr 1964\endref

\ref \no 39\by H.J. Yang \paper Transition functions and a blow-up formula
for Donaldson polynomials\paperinfo  Columbia University Thesis\yr  1992
\endref

\ref \no 40\by S.T. Yau \paper On the Ricci curvature of compact K\"ahler
manifolds and the complex Monge-Amp\`ere equation\jour Comm. Pure Appl. Math.
\vol  31 \pages  339--411 \yr 1978\endref

\endRefs
\enddocument